\documentclass[10pt]{article}
\usepackage[T1]{fontenc}
\usepackage{grffile}

\usepackage[section]{placeins}
\usepackage{chngcntr}
\usepackage{sidecap}
\usepackage[normalem]{ulem}
\usepackage[margin=.9in]{geometry}
\usepackage{amssymb, cancel, amsmath, amstext}
\usepackage[dvips]{epsfig}
\usepackage[small]{caption}
\usepackage{graphicx}
\usepackage[all]{xy}
\usepackage{mathtools}
\usepackage{hyperref}

\usepackage{caption}
\usepackage{subcaption}
\usepackage{colortbl}
\usepackage{mathrsfs}
\usepackage{tikz}
\usepackage{multirow}
\usetikzlibrary{positioning}

\usepackage[shortlabels]{enumitem}
\usepackage{wrapfig}

\usepackage{marginnote}
\usepackage{geometry}
\newtheorem{Lemma}{Lemma}[section]
\newtheorem{Theorem}[Lemma]{Theorem}
\newtheorem{Proposition}[Lemma]{Proposition}
\newtheorem{Corollary}[Lemma]{Corollary}
\newtheorem{Remark}[Lemma]{Remark}
\newtheorem{Definition}[Lemma]{Definition}

\newtheorem{Observation}[Lemma]{Observation}
\newenvironment{Proof}%
 {\begin{trivlist} \item[]{\bf Proof. }}%
 {\hspace*{\fill}$\rule{.4\baselineskip}{.4\baselineskip}$\end{trivlist}}

 {\begin{trivlist}\item[]\textbf{Acknowledgments.}}{\end{trivlist}}

\makeatletter\@addtoreset{figure}{section}\makeatother

\renewcommand{\theTheorem}{\arabic{section}.\arabic{Lemma}}

\makeatletter \@addtoreset{equation}{section} \makeatother

\newcommand{\hp}{\mathcal{H}}

\newcommand{\h}{\mathcal{H}}
\newcommand{\onelc}{(1\leadsto 0)^{(c)}}
\newcommand{\bu}{\bar{u}}
\newcommand{\tv}{\widetilde{v}}

\newcommand{\R}{\mathbb{R}}
\newcommand{\C}{\mathbb{C}}

\DeclareMathOperator{\csch}{csch}

\renewcommand{\leq}{\leqslant}
\renewcommand{\geq}{\geqslant}

\def\XXint#1#2#3{{\setbox0=\hbox{$#1{#2#3}{\int}$}
\vcenter{\hbox{$#2#3$}}\kern-.5\wd0}}

\begin{document}
\title{Horizontal patterns from finite speed directional quenching}

\author{ Rafael Monteiro\\[2ex]
\textit{\footnotesize University of Minnesota, School of Mathematics,  
206 Church St. S.E., Minneapolis, MN 55455, USA}} 

\date{\small \today} 

\maketitle

\begin{abstract}
In this paper we study the process of phase separation from directional quenching, considered as an externally triggered variation in parameters that changes the system from monostable to bistable across an interface; in our case the interface moves with speed $c$ in such a way that the bistable region grows. According to results from  \cite{Monteiro_Scheel-contact_angle,Monteiro_Scheel}, several patterns exist when $c \gtrsim0$, and here we investigate their persistence for finite $c>0$, clarifying the pattern selection mechanism related to the speed $c$ of the quenching front.  
\vspace{\baselineskip}

\end{abstract}

\section{Introduction}

\begin{figure}[htb]
    \centering
    \begin{subfigure}[b]{0.45\textwidth}
          \centering
        \includegraphics[height = 0.4in]{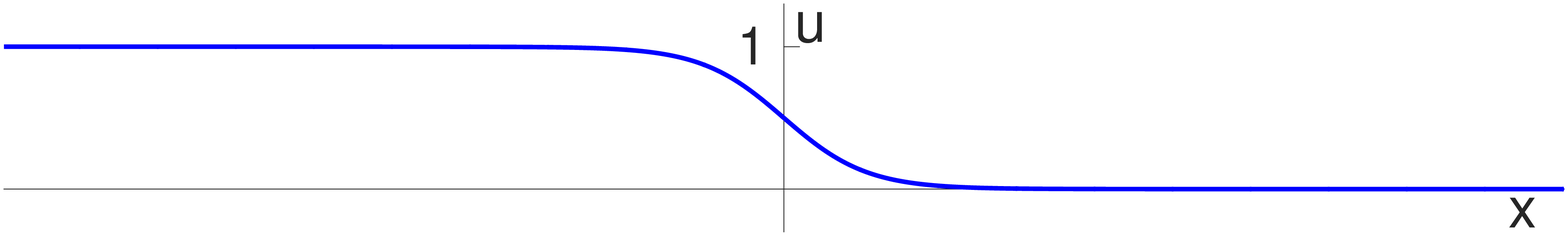}
    \end{subfigure}
    ~ 
    \begin{subfigure}[b]{0.45\textwidth}
    \centering
    \includegraphics[width=\textwidth,height=0.4in]{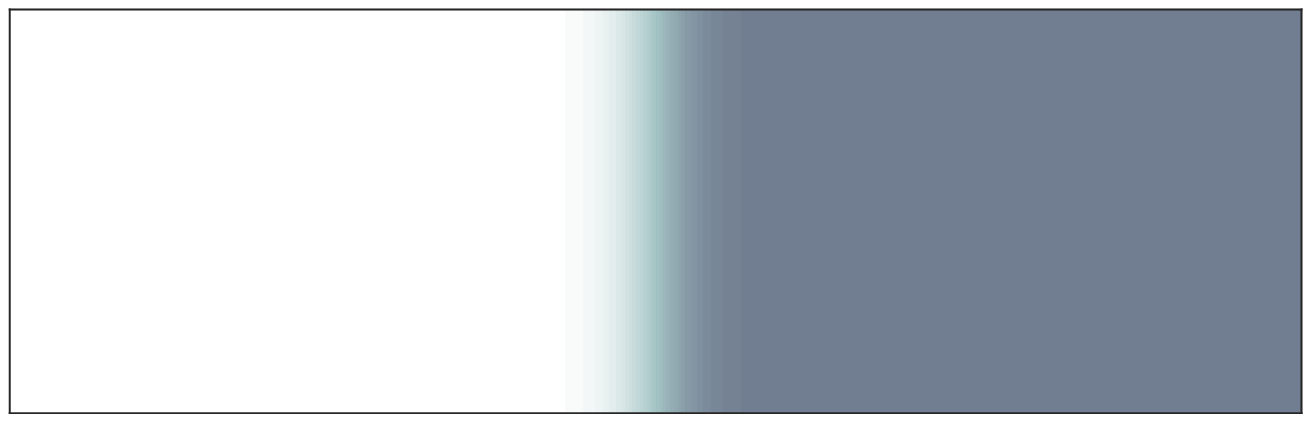}
        \end{subfigure}
    \caption{Sketches of solutions for pure phase selection $1 \leadsto 0$; solution $\theta^{(c)}(x)$ (left) and contour plot for $(x,y)\in\R^2$ (right).\label{Figure_1-1_leadsto_0}}
\end{figure}

\begin{figure}[htb]
    \centering
    \begin{subfigure}[b]{0.45\textwidth}
          \centering
    \includegraphics[width=\textwidth,height=0.4in]{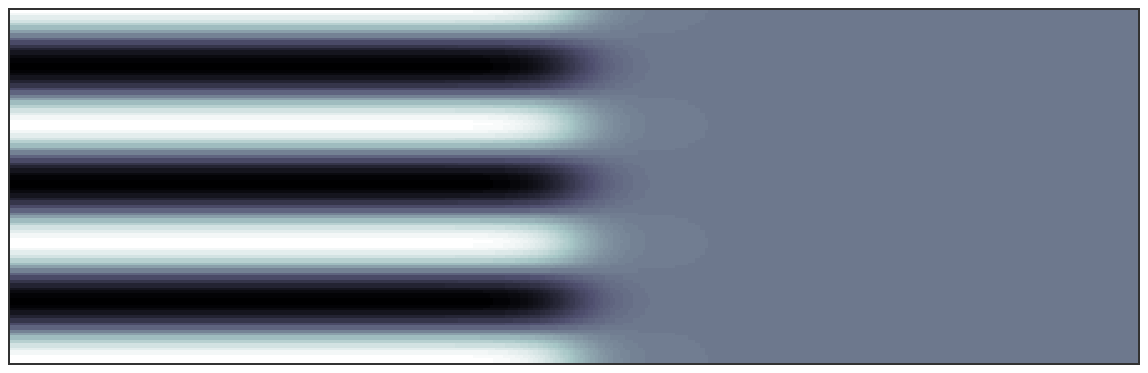}    
              \end{subfigure}
    ~ 
    \begin{subfigure}[b]{0.45\textwidth}
    \centering
   \includegraphics[width=\textwidth,height=0.4in]{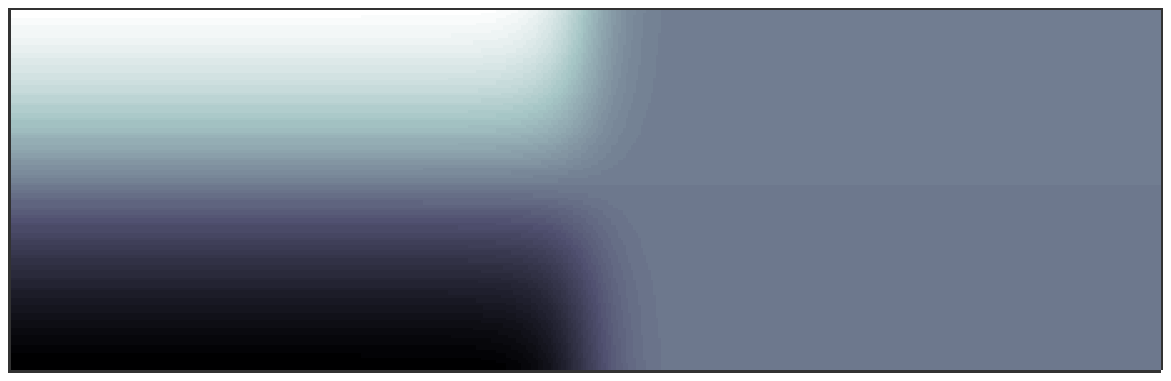} 
       \end{subfigure}
   \caption{Sketches of solutions for horizontal patterns;  $\hp_{\kappa}$ pattern (left) and $\mathcal{H}_{\infty}$ pattern (right). \label{Figure_1-horizontal}}
\end{figure}

In the theory of reaction diffusion, the interplay between stable and unstable mechanisms can give rise to spatial patterns, i.e., stationary non-homogeneous structures. In the presence of controllable external parameters  the existence and persistence of these patterns are worth to investigate, both for their mathematical interest and industrial applications; see \cite{Langmuir} and the survey \cite{Foard}. In this paper we are interested in a \textit{directional quenching} scenario, where a planar interface (also called \textit{the quenching front}) moves with constant speed $c$, across which a phase separation process takes place: ahead of the interface the system is monostable, while in its wake it is bistable. We study this phenomenon in  the scalar model
 \begin{eqnarray}\label{directional_quenching_eq:not_moving_frame}
\partial_tu  = \Delta_{(\xi,y)} u + \mu(\xi - ct) u  - u^3,
\end{eqnarray}
where $\Delta_{(\xi,y)}:= \partial_{\xi}^2 + \partial_y^2$ and  $\mu(\xi - ct \gtrless 0) = \mp 1$. Equation \eqref{directional_quenching_eq:not_moving_frame} is a particular case of the \textit{Allen-Cahn} model, which describe the behavior of a heterogeneous, binary mixture: the unknown $u(\xi,y; t)$ denotes the relative concentration of one of the two metallic components of the alloy at time $t\in \R^+ :=[0,\infty)$ and point $(\xi,y)\in \R^2$.  The stationary problem in a moving frame $(x,t) = (\xi - ct, t)$ is written 
 \begin{eqnarray}\label{Main_equation}
 - c\partial_x u = \Delta_{(x,y)} u + \mu(x) u  - u^3 .
\end{eqnarray}
The most physically relevant scenario to be considered consists of the case $c\geq 0$. According to results in \cite{Monteiro_Scheel}, whenever the speed $c$ of the quenching front is small the equation \eqref{Main_equation} admits a rich family of patterns , as we now briefly describe. 
\paragraph{Pure phase selection: \texorpdfstring{\textbf{$\onelc$}}{1->0} fronts.}
This is the simplest, one dimensional case, when \eqref{Main_equation} reads
\begin{align}\label{1_leadsto_0:1D}
- c\partial_x u(x) = \partial_x^2 u(x) + \mu(x)u(x)  - u(x)^3, \qquad x\in \R.
\end{align}
The quenching trigger generates a pattern $\theta^{(c)}(\cdot)$ solving \eqref{1_leadsto_0:1D} and satisfying spatial asymptotic conditions $\displaystyle{\lim_{x\to -\infty}\theta^{(c)}(x)= 1}$  and  $\displaystyle{\lim_{x\to +\infty}\theta^{(c)}(x)= 0}$ in the wake and ahead of the quenching front, respectively; see Figure \ref{Figure_1-1_leadsto_0}. 
\paragraph{Horizontal patterns: \texorpdfstring{$\mathcal{H}_{\kappa}$, $\pi <\kappa \leq \infty$}{H_kappa, pi <k <=infty}.}
In this scenario the patterns  sought are truly two-dimensional; furthermore, we can take advantage of the odd nonlinearity to reduce the study of \eqref{Main_equation} to a  problem in  $\R\times[0,\kappa]$. Whenever $\kappa <\infty$ the solution $\displaystyle{\Xi_{\kappa}^{(c)}(\cdot,\cdot)}$ has boundary conditions $\displaystyle{\Xi_{\kappa}^{(c)}(x,y)\Big|_{y=0,\kappa}=0}$ and  spatial asymptotic conditions  
\begin{align}\label{H_kappa:asymptotics}
 \lim_{x\to +\infty}\Xi_{\kappa}^{(c)}(x,y)= 0, \quad \mbox{and} \quad \lim_{x\to -\infty}|\Xi_{\kappa}^{(c)}(x,y)-\bar{u}(y;\kappa)|= 0,
\end{align}
where  $\bar{u}(y;\kappa)$ is a parametrized family of periodic solutions to
\begin{equation}\label{periodic_solutions}
\partial_y^2\bar{u}+\bar{u}-\bar{u}^3=0, \qquad \bar{u}(y;\kappa)=-\bar{u}(y+\kappa;\kappa)=-\bar{u}(-y;\kappa)\not\equiv 0,\qquad \mbox{for}\quad y\in\R.
\end{equation}
with half-periods  $\pi<\kappa<\infty$, and normalization $\partial_y\bar{u}(0;\kappa)>0$, $\bar{u}(0;\kappa) = 0$. In the limiting case $\kappa = \infty$ a single interface is created, describing what we call $\hp_{\infty}$-pattern (see Fig. \ref{Figure_1-horizontal}); asymptotically we have
\begin{align}\label{H_infty:asymptotics}
 \lim_{y\to +\infty}\Xi_{\infty}^{(c)}(x,y)= \theta^{(c)}(x), \quad \lim_{x\to +\infty}\Xi_{\infty}^{(c)}(x,y)= 0, \quad \mbox{and} \quad \lim_{x\to -\infty}\Xi_{\infty}^{(c)}(x,y)=u(y;\infty),
\end{align}
where $u(\cdot;\infty):=\tanh\left(\frac{\cdot}{\sqrt{2}}\right)$ and $\theta^{(c)}(\cdot)$ is a $\onelc$-front.

We summarize  below the properties of these solutions in the regime $c \gtrsim 0$ as given in \cite{Monteiro_Scheel}. Their approach is based on a continuation argument from the case in which quenching front has zero speed ($c=0$), somehow explaining the nature of the smallness on $c$ in their results.%
\begin{Proposition}[{\cite[Theorems 1.1 and 1.5]{Monteiro_Scheel}}; \textbf{$\onelc$ problem}]\label{DQ:c_zero:1_leadsto_0} For any $c \geq 0$ sufficiently small there exists a function $\theta^{(c)}(\cdot) \in \mathscr{C}^{(1,\alpha)}(\R; [0,1])$, $\forall \alpha \in [0,1)$, that solves the $\onelc$ problem, i.e., $\theta^{(c)}$ solves the boundary value problem  
\begin{eqnarray*}
\left\{\begin{array}{c}
  \partial_x^2\theta^{(c)}(x) + c\partial_x\theta^{(c)}(x) + \mu(x) \theta^{(c)}(x) - \left(\theta^{(c)}(x)\right)^3 = 0 \quad \mbox{in the sense of  distributions} \\
 0< \theta^{(c)}(x) <1, \quad \theta^{(c)}(-\infty)=1, \quad \theta^{(c)}(+\infty)=0, 
\end{array}\right.
\end{eqnarray*}
where the boundary conditions are attained in the limit sense. Furthermore, the mapping $x \mapsto \theta^{(c)}(x)$ is non-increasing. The solutions found for $c=0$ can be continued smoothly to $c>0$ for sufficiently small $c$. More precisely, there exist families of solutions  $\theta^{(c)}(x)$ to \eqref{1_leadsto_0:1D} for $0<c<\delta_1$, satisfying the same limiting conditions as the solutions at $c=0$ for $x\to\pm\infty$. Moreover, the solutions depend smoothly on $c$, uniformly in $x$.
 \end{Proposition}
\begin{Proposition}[{\cite[Theorem 1 and Proposition 1.4]{Monteiro_Scheel}}; $\hp_{\kappa}$ problem, $\pi < \kappa \leq \infty$]\label{DQ:c_zero:H_infty} For any $c \geq 0$ sufficiently small equation \eqref{Main_equation} admits a family of solutions (in the sense of distributions) $\Xi_{\kappa}^{(c)}(\cdot,\cdot)\in \mathscr{C}^{(1,\alpha)}(\R^2; \R)$, $\forall \alpha \in [0,1)$, $\kappa\in (\pi,\infty]$. In the case $\pi < \kappa <\infty$ for the solution $\Xi_{\kappa}^{(c)}(\cdot, \cdot)$ has the symmetries 
$\Xi_{\kappa}^{(c)}(x,y)=-\Xi_{\kappa}^{(c)}(x,-y)=-\Xi_{\kappa}^{(c)}(x,y+\kappa)=-\Xi_{\kappa}^{(c)}(x,y+\kappa)$ and satisfies  the asymptotic spatial conditions \eqref{H_kappa:asymptotics}. Moreover, the convergence is exponential, uniformly in $y$. On the other hand, whenever $\kappa = \infty$ the solution $\Xi_{\infty}^{(c)}(\cdot, \cdot)$ has the symmetries $\Xi_{\infty}^{(c)}(x,y)=-\Xi_{\infty}^{(c)}(x,-y)$ and satisfies the asymptotic spatial conditions \eqref{H_infty:asymptotics}, where convergence is exponential and  uniform.
\end{Proposition}
In this paper we give a deeper understanding of the range of validity of these continuations in $c$. 
\begin{Remark}
Besides the patterns described above, oblique and vertical structures with respect to the quenching front were also studied in \cite{Monteiro_Scheel}, where they were shown to not exist as solutions to \eqref{Main_equation} when $c> 0$. Therefore, the only patterns of relevance for us are the ones described above. Throughout this paper we add sub and superscripts to the patterns found in \cite{Monteiro_Scheel} under the inconvenience of disagreeing with that paper's notation; this is done because the classification of patterns has become more involved and richer. In this way our notation highlights the dependence of the solutions on the quenching speed $c$ and on the  $\kappa$-periodicity in the y-direction as  $x \to -\infty$ (in the multidimensional case; see Fig. \ref{Figure_1-horizontal}). 
\end{Remark}

\subsection{Main results}

Initially we study the problem in one dimensional setting; although less physically relevant, it stands as one of the cornerstone of the construction of the 2D patterns -- $\mathcal{H}_{\kappa}$ and $\mathcal{H}_{\infty}$ (see for instance, \cite[Sec. 2 and 3]{Monteiro_Scheel}).

\begin{Theorem}[$\onelc$-patterns]\label{theorem:1D_result} For any fixed $c \in [0,2)$ there exists a unique solution $\theta^{(c)}(\cdot) \in \mathscr{C}^{(1,\alpha)}(\R; (0,1))$, $\forall \alpha \in [0,1)$ satisfying
\begin{eqnarray}\label{theorem:1D_result-equation}
\left\{\begin{array}{c}
  \partial_x^2\theta^{(c)}(x) + c \partial_x\theta^{(c)}(x) + \mu(x)\theta^{(c)}(x) - (\theta^{(c)})^3(x) = 0 \quad \mbox{in the sense of  distributions} \\
 0< \theta^{(c)}(x) <1, \displaystyle{\quad \lim_{x\to -\infty}\theta^{(c)}(x)=1, \quad \lim_{x\to \infty}\theta^{(c)}(x)=0}. 
\end{array}\right.
\end{eqnarray}
The convergence takes place at exponential rate. Furthermore, the solution solution $\theta^{(c)}(\cdot)$ is strictly positive and  satisfies $(\theta^{(c)})'(x) <0$ for all $x \in \R$. No solution to \eqref{theorem:1D_result-equation} exists when $c >2$.
\end{Theorem}
In the higher dimensional setting, both  the magnitude of the speed $c$ of the quenching front and the y-period $\kappa$ of the end-state $\bar{u}(\cdot; \kappa)$ play important roles in the analysis.


\begin{SCfigure}[][htb]
  \centering
\includegraphics[width=0.25\textwidth]{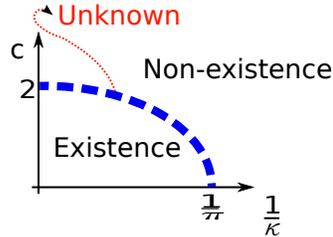}
\caption{Existence diagram  for parameters $c\geq0$ (speed of the front) and $\kappa  > \pi$ (y-periodicity of the patterns); the dashed curve represents the critical case $\mathcal{P}(c,\kappa) =1$ (see Theorem \ref{theorem:2D_result}). \label{Fig:c_versus_kappa}}
\end{SCfigure}
\begin{Theorem}[$\mathcal{H}_{\kappa}$, $\mathcal{H}_\infty$ patterns]\label{theorem:2D_result}
Let $\pi <\kappa \leq\infty$ be a fixed number. 
Define the quantity
\begin{align}\label{critical_quantity}
 \mathcal{P}(c;\kappa) := \left\{\begin{array}{ccc}
                           \frac{c^2}{4} + \frac{\pi^2}{\kappa^2},&\quad  \mbox{for} \quad \pi<\kappa <\infty; \\
                           \frac{c^2}{4},&\quad  \mbox{for} \quad \kappa = \infty.
                          \end{array}\right.
\end{align}
\begin{enumerate}[label=(\roman*), ref=\theTheorem(\roman*)]
\hangindent\leftmargin
 \item \textbf{(Existence)} \label{theorem:2D_result:existence} The solutions $\Xi_{\kappa}^{(0)}(\cdot, \cdot)$ defined in Prop. \ref{DQ:c_zero:H_infty} when  $c=0$ can be continued smoothly in $c>0$ within the range $\mathcal{P}(c;\kappa) < 1$ to  solutions  $\Xi_{\kappa}^{(c)}(\cdot,\cdot)$  solving \eqref{Main_equation} in the sense of distributions and satisfying the asymptotic spatial condition \eqref{H_kappa:asymptotics} (resp. \eqref{H_infty:asymptotics}) when $\pi <\kappa <\infty$ (resp., when $\kappa = \infty$). The convergence to their spatial asymptotic states takes place at exponential rate, uniformly. Furthermore, for any $\kappa > \pi$ the mapping $x\mapsto \Xi_{\kappa}^{(c)}(x,y)$ is non-increasing for any fixed $y \in[0,\kappa]$ and $0 <\Xi_{\kappa}^{(c)}(x,y) <\bar{u}(y;\kappa)$ in $(x,y) \in \R\times(0,\kappa)$, where $\bar{u}(\cdot,\cdot)$ is given by  \ref{periodic_solutions}. The solutions have the symmetries $\Xi_{\kappa}^{(c)}(x,y)=-\Xi_{\kappa}^{(c)}(x,-y)=-\Xi_{\kappa}^{(c)}(x,y+\kappa)=-\Xi_{\kappa}^{(c)}(x,y+\kappa)$ (resp., $ \Xi_{\infty}^{(c)}(x,y)=-\Xi_{\infty}^{(c)}(x,-y)$). 

 \item \textbf{(Nonexistence)}\label{theorem:2D_result:nonexistence} Whenever $\kappa \in (\pi, \infty)$ (resp., $\kappa =\infty$) and $\mathcal{P}(c;\kappa) > 1$ no solution to \eqref{Main_equation} satisfying $0 <\Xi_{\kappa}^{(c)}(x,y) <\bar{u}(y;\kappa)$ in $(x,y) \in \R\times(0,\kappa)$ and the asymptotic spatial condition \eqref{H_kappa:asymptotics} (resp. \eqref{H_infty:asymptotics}) can be found. 
\end{enumerate}
\end{Theorem}
One can see from the previous result that whenever $c \geq 0$ the region $\mathcal{P}(c;\kappa) <1$ (resp. $\mathcal{P}(c;\kappa) >1$) corresponds  to $\displaystyle{c < \sqrt{1 - \frac{\pi^2}{\kappa^2}}}$ (resp., $\displaystyle{c > \sqrt{1 - \frac{\pi^2}{\kappa^2}}}$) , namely, the linear spreading speed obtained from the linearization of \eqref{Main_equation} about $u\equiv  0$ on the region $x\leq 0$. Overall, we point out that the dependence of the critical spreading speed curve on the parameter $\kappa$ is a true manifestation of the multidimensionality of the $\mathcal{H}_{\kappa}$-patterns; the quantity $\mathcal{P}(\cdot;\cdot)$ describes the maximal speed of spreading of the bistable region and its dependence on the  y-period of the pattern away from the quenching interface. 

\begin{Remark}[Uniqueness results]
 It is worth to point out that the uniqueness result in Theorems \ref{theorem:1D_result} and \ref{theorem:2D_result} allow us to compare the solutions constructed in \cite{Monteiro_Scheel} using perturbation methods for $c \gtrsim 0$ with those obtained here.
\end{Remark}

\paragraph{Critical cases; \texorpdfstring{ $\mathcal{P}(c;\kappa) = 1$.}{S=1}} To the best of the author's knowledge, the cases $c = 2$ (for the $\onelc$-problem), and $\displaystyle{\mathcal{P}(c;\kappa) =1}$ (for the $\hp_{\kappa}$ and $\hp_{\infty}$-problems) are open. See Sec. \ref{Discussion}.

A summary of our results is given in the Table \ref{table-1}.
%

  \begin{table}[!htbp]
\begin{minipage}{0.48\textwidth}
 \begin{center}
\begin{tabular}{ccc}
\hline \\
\multicolumn{3}{c}{\begin{minipage}{3cm}
\begin{center}
$\onelc$ problem
\end{center}
\end{minipage} }\\
[1em]
\hline\\

$0 \leq c< 2$ & $c = 2$ & $c >2$ \\
 [1em]
\hline\\
 \begin{minipage}{2cm}
\begin{center}
\textbf{Yes}\\
Thm. \ref{theorem:1D_result}
\end{center}
\end{minipage}
& 
\begin{minipage}{2cm}
\begin{center}
\textbf{Not known}
\end{center}
\end{minipage}
&
\begin{minipage}{2cm}
\begin{center}
\textbf{No}\\
Thm. \ref{theorem:1D_result}
\end{center}
\end{minipage}
\\
\\
 \hline 
\end{tabular}
\end{center}
\end{minipage}
\hfill 
\begin{minipage}{0.48\textwidth}
\begin{center}
\begin{tabular}{ccc}
  \hline \\
  \multicolumn{3}{c}{\begin{minipage}{4cm}
\begin{center}
$\h_{\kappa}$ problem ($\pi <\kappa \leq\infty$)
\end{center}
\end{minipage}}\\
[1em]

\hline\\
$\mathcal{P}(c;\kappa)< 1$ & $\mathcal{P}(c;\kappa) = 1$ & $\mathcal{P}(c;\kappa) >1$ \\
 [1em]
\hline\\
 \begin{minipage}{2cm}
\begin{center}
\textbf{Yes}\\
Thm. \ref{theorem:2D_result:existence}
\end{center}
\end{minipage}
& 
\begin{minipage}{2cm}
\begin{center}
\textbf{Not known}
\end{center}
\end{minipage}
&
\begin{minipage}{2cm}
\begin{center}
\textbf{No}\\
Thm. \ref{theorem:2D_result:nonexistence}
\end{center}
\end{minipage}
\\
\\
 \hline 
\end{tabular}
\end{center}
\end{minipage}
\caption{Existence tables for the patterns $\onelc$ and $\mathcal{H}_{\kappa}$ ($\pi < \kappa \leq \infty$), for critical quantity $\mathcal{P}(c;\kappa)$, defined in \eqref{critical_quantity}.\label{table-1}}
 \end{table}

\paragraph{Single interfaces with contact angle.} The  $\hp_{\infty}$-patterns obtained in Theorem \ref{theorem:2D_result} are odd functions with respect to the $y$ variable, and as so they satisfy $\Xi_{\infty}^{(c)}(x,0) = 0$. Thus, the patterns present a nodal set at the negative part of the $x$-axis that forms a  right angle $\left(\varphi = \frac
{\pi}{2}\right)$ with respect to the quenching front located at a line $x- ct=0$ (parallel to the $y$ axis). This observation motivates the question of how extra terms added to equation \eqref{Main_equation} could deform this nodal set. More precisely, we consider the equation 
\begin{equation}\label{unbalanced_equation}
\Delta u +c_x\partial_x u + c_y\partial_y u +  \mu(x) u  - u^3 + \alpha g(x,u)=0 ,\qquad (x,y)\in\R^2.
\end{equation}
where $\mu(x) = \pm 1$ when $x \lessgtr 0$
and $\displaystyle{g(x,u)=\left\{\begin{array}{ll}g_\mathrm{l}(u),&x<0\\g_\mathrm{r}(u),&x>0\end{array}\right.
},$ for $\displaystyle{g_{\{l,r\}}(\cdot) \in \mathscr{C}^{\infty}(\R;\R)}$. Notice that $c_y = 0$, $\alpha =0$, the function $u(\cdot,\cdot) = \Xi_{\infty}^{(c_x)}(\cdot,\cdot)$ solves \eqref{unbalanced_equation}.

With regards to \eqref{unbalanced_equation}, we remark that two mechanisms are in play:  the growth of the bistable region and new, ``unbalancing terms'', that break the odd symmetry of the solutions. 

\begin{Definition}[Contact angle]\label{definition:contact_angle}
%
We say \eqref{unbalanced_equation} possesses a solution $u$ with contact angle $\varphi_*$ if $u$ possesses the limits 
\begin{equation}\label{unbalanced_limit}
\lim_{x\to +\infty} u(x,y)=0, \qquad \lim_{x\to -\infty} u(x,\cot(\varphi)x)=\left\{\begin{array}{ll}
u_+,& \varphi>\varphi_*\\
u_-,& \varphi<\varphi_*
\end{array}\right. ,
\end{equation}
for all $0<\varphi<\pi$.  
\end{Definition}
For instance, when $\alpha = 0$ and $c_y=0$ the function $\Xi_{\infty}^{(c_x)}(\cdot,\cdot)$ solves \eqref{unbalanced_equation} satisfying the limit \eqref{unbalanced_limit} for $\varphi^{*} = \frac{\pi}{2}$ and $u_{\pm} = \pm1$.  It was shown in \cite{Monteiro_Scheel-contact_angle} shows that  for any $c_x\gtrsim0$ fixed, the function $\Xi_{\infty}^{(c=c_x)}(\cdot,\cdot)$ can be continued in $\varphi^{*}$  as a solution to \eqref{unbalanced_equation} for all $\left|\varphi^{*} - \frac{\pi}{2}\right|$ sufficiently small. One of the most important properties used  in their proof concerns to the  strict monotonicity of the mapping $y \mapsto \Xi_{\infty}^{(c_x)}(x,y)$ for any $x\in \R$, namely,
\begin{equation}\label{unbalanced_equation:monotonicity}
 \partial_y\Xi_{\infty}^{(c_x)}(x,y) > 0, \quad \mbox{for} \quad (x,y) \in \R^2.
\end{equation}
where $\Xi_{\infty}^{(c_x)}(\cdot,\cdot)$ is given in Theorem \ref{theorem:2D_result:existence}. An inconvenient in their analysis is the fact that the patterns $\Xi_{\infty}^{(c_x)}$ for $c_x>0$ were obtained through perturbation methods (in \cite[\S 5]{Monteiro_Scheel}), hence some qualitative information on the patterns are not immediately available. Nevertheless, the authors managed to prove \eqref{unbalanced_equation:monotonicity} for $c_x \gtrsim0$ sufficiently small (see \cite[Prop. 4.1]{Monteiro_Scheel-contact_angle}).
Our construction readily gives the validity of \eqref{unbalanced_equation:monotonicity} for all $0<c_x<2$, thus we can  make use of the analysis in \cite{Monteiro_Scheel-contact_angle} to conclude the following result:
\begin{Corollary}[Unbalanced patterns]
 For $0<c_x < 2$, there exists $\alpha_0(c_x)>0$ such that for all $|\alpha|<\alpha_0(c_x)$ there exist a speed $c_y(\alpha)$ and a solution  $u(x,y;\alpha)$ to \eqref{unbalanced_equation} with contact angle $\varphi(\alpha)$. We have that $u(x,y;\frac{\pi}{2}) = \Xi^{(c)}(\cdot, \cdot; \infty)$. Moreover, $c_y(\alpha)$ and $\varphi(\alpha)$ are smooth with $\varphi(0)=\pi/2$, $c_y(0)=0$, and $u(x,y;\alpha)$ is smooth in $\alpha$  in a locally uniform topology, that is, considering the restriction $u|_{B_R(0)}$ to an arbitrary large ball.

\end{Corollary}

\begin{figure}[htb]
\begin{center}
\includegraphics[width=0.3\textwidth]{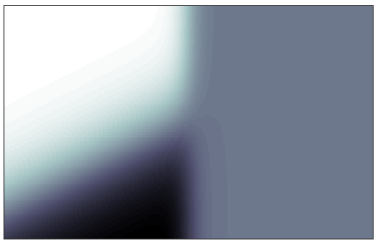}
\end{center}
\caption{Sketch of an unbalanced pattern with a contact angle; see Def. \ref{definition:contact_angle} or \cite{Monteiro_Scheel-contact_angle}. \label{Fig:contact_angle}}
\end{figure}

\subsection{Outline}

In Section \ref{Sec:H_1d_and_1_0_problem:truncation_approach} we focus on the $\onelc$ problem, where we prove Theorem \ref{theorem:1D_result}. Section \ref{Sec:2D_result:k_finite} is devoted to Theorem \ref{theorem:2D_result} and $\hp_{\kappa}$ patterns ($\pi <\kappa <\infty$), while the study of the $\hp_{\infty}$-pattern is left to  Section \ref{Sec:2D_result:k_infinite}.  A brief discussion and further extensions brings the paper to an end in Section \ref{Discussion}. 

\paragraph{Notation.} In this paper we write $\mathscr{C}^k(X;Y)$, $\mathscr{C}_0^k(X;Y)$ and $\mathscr{C}^{(k,\alpha)}(X;Y)$ denote respectively, the space of $k$-times continuously differentiable functions, the space of $k$ times continuously differentiable functions with compact support in $X$, the space of $(k, \alpha)$-H\"older continuously differentiable functions from $X$ to $Y$. We denote the Sobolev spaces over an open set $\Omega$  by $H^k(\Omega)$. The inner product of elements in a Hilbert space $\mathcal{H}$ is written as $\langle, \rangle_{\mathcal{H}}$. Norms on a Banach space $\mathcal{B}$ are denoted as $||\cdot||_{\mathcal{B}}$.
For a given operator $\mathscr{L}: \mathcal{D}(\mathscr{L}) \subset X \to Y$ we write $\mathrm{Ker}\left( \mathscr{L} \right) :=\{ u \in \mathcal{D}(\mathscr{L})| \mathscr{L}u=0 \}$ and $\mathrm{Rg}\left( 
\mathscr{L} \right) :=\{ f \in Y |\exists u \in \mathcal{D}(\mathscr{L}), Lu =f \}.$ A distribution $T \in \mathcal{
D}'(\Omega)$ satisfies  $T\geq 0$ in the sense of distributions if $T(\phi) \geq 0$ for any $\phi(\cdot) \in \mathcal{C}_0^{\infty}(\Omega; [0, \infty)).$
We define a $\mathscr{C}^{\infty}(\R;[0,1])$ partition of unity $\{\chi^{\pm}(\cdot)\}$ of $\R$, of the form 
\begin{align}\label{partition_of_unity}
\chi^{-}(x) + \chi^{+}(x)  =1,\quad  \mbox{where} \quad \chi^{-}(x) =1 \quad \mbox{for}\quad  x\leq -2,\quad \mbox{and} \quad \chi^{-}(x) =0\quad \mbox{for} \quad  x\geq -1.  
\end{align}
Last, we denote the Implicit Function Theorem by IFT. 
\paragraph{Acknowledgments.} The author is grateful to Prof. Arnd Scheel for many interesting discussions and insights throughout the writing of this paper. Many thanks also go  to Prof. Yasumasa Nishiura and Prof. Natsuhiko Yoshinaga for sharing their perspective on this work.  The author acknowledges partial support through  NSF grants DMS-1612441 and DMS-1311740.
 \section{One dimensional directional quenching: \texorpdfstring{$(1\leadsto 0)^c$}{1-0} problem, \texorpdfstring{$c >0$}{c>0}}\label{Sec:H_1d_and_1_0_problem:truncation_approach}

The construction of the patterns $\onelc$ follows ideas from  \cite{Monteiro_Scheel} and \cite{kolli2003approximation}: initially we solve a family of similar problems in truncated, bounded intervals; later, as we enlarge these intervals and exhaust $\R$, we show that these functions  converge to  a solution of problem \eqref{1_leadsto_0:1D}. We begin by setting up a \textit{truncated $(1\leadsto 0)^c$ problem}:
\begin{eqnarray}\label{H_1D_problem:truncated}
\left\{\begin{array}{cc}
 \partial_x^2u(x) + c\partial_xu(x) + \mu(x) u(x)  - u^3(x) &= 0, \\
 u(-M)= 1,  u(L)=0, &
\end{array}\right. 
\end{eqnarray}
for $0<M,L$, with continuity of $u$ and $u_x$ at $x=0$. It is shown the existence of a unique solution $\theta^{(c)}_{(-M,L)}$; later on the section we  let $M\to\infty$ and, subsequently, $L\to\infty$.  Roughly speaking, the $\onelc$ front $\theta^{(c)}(\cdot)$  will be given by
\begin{align}\label{minimax}
\displaystyle{\theta^{(c)}(\cdot) = \sup_{L>0} \left\{\inf_{M>0}\theta^{(c)}_{(-M,L)}(\cdot) \right\}}. 
\end{align}
%
The qualitative properties of the function $\theta^{(c)}(\cdot)$ are proved in this section, where we also show that $\theta^{(c)}(x)$ converges to $1$ and $0$ as $x\to\infty$ and $x\to-\infty$, respectively. We finalize with a proof of Theorem \ref{theorem:1D_result}.  A substantial part of the techniques and proofs are similar to those in \cite{Monteiro_Scheel}; whenever possible we skip details and refer to that article for full proofs. 
\subsection{The \texorpdfstring{$\onelc$}{1->0c} truncated problem} \label{s:21}
\begin{Lemma}[Existence and uniqueness]\label{Existence_uniqueness_truncated:1_leadsto_0} The truncated problem \eqref{H_1D_problem:truncated} has a unique solution $\theta_{(-M,L)}^{(c)}(\cdot) \in \mathscr{C}^{(1,\alpha)}([-M,L];[0,1])$,  $\forall \alpha \in [0,1)$. Furthermore, $\theta_{(-M,L)}^{(c)}(x) >0$ whenever $x \in (-M,L)$.\end{Lemma}
\begin{Proof}
To prove existence of a solution we define an iterative scheme,
\begin{eqnarray*}
\left\{\begin{array}{c}
  -  U_{n+1}^{''}(x) - cU_{n+1}'(x) +5U_{n+1}(x)  =  \left(5 +\mu(x) \right) U_n(x) -U_n^3(x), \quad \mbox{in} \, (-M,L) \\
 U_{n+1}(-M)=1, \quad U_{n+1}(L)=0,
\end{array}\right.
\end{eqnarray*}
where we write $(\cdot)' = \partial_x(\cdot)$. Following the reasoning in \cite[Sec. 2.2]{Monteiro_Scheel}, it is shown that $\{U_n\}_{n\in \mathbb{N}}$ is pointwisely decreasing and so that $\displaystyle{\theta_{(-M,L)}^{(c)}(x) := \inf_{n \in \mathbb{N}}U_n(x)}$ is a well defined element of $ \displaystyle{\mathscr{C}^{(1,\alpha)}([-M,L];[0,1])}$,  $\forall \alpha \in [0,1)$. Furthermore, $\displaystyle{\theta_{(-M,L)}^{(c)}(x) >0}$ for $x \in (-M,L)$.  The uniqueness proof is a bit different due to the  transport term $c \partial_x$ and we give it here for completeness: assume the existence of  two solutions, $\theta(\cdot)$, $\widetilde{\theta}(\cdot)$ so that $\theta (\cdot)\not\equiv \widetilde{\theta}(\cdot)$. Define the set  $\mathscr{I} = \{x \in [-M,L] \, |\,  \theta(x) \neq \widetilde{\theta}(x) \};$ this set is open due to continuity of $\theta$, $\widetilde{\theta}$. As an open subset of the real line, we can assume without loss of generality that  $\mathscr{I} = (a,b)$ where $\theta(x)> \widetilde{\theta}(x)$, $x \in (a,b)$, and  $\theta(x)= \widetilde{\theta}(x)$, $x\in\{a,b\}$. Now, since both $\theta$ and $\widetilde{\theta}$ are solutions, we can integrate against test functions\footnote{Recall that distributions of finite order (say, order $k$) can be extended to the space of $\mathscr{C}_0^k$ functions (cf. \cite[\S 2]{Hormander}; see also \cite[Lem. 2.2]{Monteiro_Scheel}).} $e^{c\,x}\theta(x)$ and $e^{c\,x}\widetilde{\theta}(x)$ on the interval $(a,b)$:
 \begin{align*}
 \int_a^b e^{c\, x}(-\theta''\widetilde{\theta} +\theta \widetilde{\theta}'')(x) dx - c\int_a^b e^{c\, x}(\theta' \widetilde{\theta} - \theta \widetilde{\theta}')(x)dx = \int_a^be^{c\, x}(\theta^2 -  \widetilde{\theta}^2) (x)\theta(x) \widetilde{\theta}(x)dx, 
 \end{align*}
 where $(\cdot)'$ denotes $\partial_x$. Integration by parts gives
 \begin{align*}
 e^{c\, x}(-\theta'\widetilde{\theta} +\theta \widetilde{\theta}')(x)\Big|_a^b = \int_a^be^{c\, x}(\theta^2 -  \widetilde{\theta}^2)(x) \theta(x) \widetilde{\theta}(x)dx. 
 \end{align*}
The term on the left hand side is non-positive, since $\theta > \widetilde{\theta}$ in $(a,b)$, $\theta(x) = \widetilde{\theta}(x)$ for $x \in \{a,b\}$. On the right hand side,  the term $\theta\widetilde{\theta}$ is strictly positive, thanks to the strict positivity of solutions in $(-M,L)$.  Using that $\theta > \widetilde{\theta}$ in $(a,b)$ we conclude that the integral on the right hand side is positive. This contradiction proves the result. 
\end{Proof}
In order to compare the families of solutions as $M,L$ vary, we construct trivial extensions of functions $u$ defined on an interval $(-M,L)$ given by the  operator $\mathscr{E}$, 
\begin{eqnarray*}
 \mathscr{E}\left[u\right](x) = \left\{ \begin{array}{ccc}
                        u(x), & \quad & \mbox{for} \quad x \in (-M, L)\\
                        1, &  \quad & \mbox{for}\quad  x \leq -M\\
                         0, &  \quad & \mbox{for} \quad  x \geq L.
                       \end{array} \right.
\end{eqnarray*}
By construction, $\mathscr{E}\left[\theta_{(-M,L)}^{(c)}(x)\right] $ is a continuous function, and this extension is one of the main tools to make \eqref{minimax} meaningful and rigorous. The proof of the next result  follows the results in \cite[\S 2]{Monteiro_Scheel}.

\begin{Lemma}[Properties of solutions to the truncated problem]\label{Thm:monotonicity_properties_1_to_0} The following properties of $\mathscr{E}[\theta_{(-M,L)}^{(c)}](\cdot)$ hold.
\begin{enumerate}[label=(\roman*), ref=\theTheorem(\roman*)]
\hangindent\leftmargin
 \item\emph{(Monotonicity of $\mathscr{E}$)}\label{Thm:monotonicity_properties_1_to_0:monotonicity_extension_operator}	
 We have $0 \leq \mathscr{E}\left[\theta_{(-M,L)}^{(c)}\right](\cdot) \leq 1$. Furthermore, for $w$ defined on a subset $A$, $ (-M,L) \subset A \subset (-\infty,L)$ with  $0 \leq w(\cdot) \leq 1$ and $ 0 \leq w(\cdot) \leq \theta_{(-M,L)}^{(c)}(\cdot)$ in $(-M,L)$, we have
$0 \leq \mathscr{E}\left[w\right](\cdot) \leq \mathscr{E}\left[\theta_{(-M,L)}^{(c)}\right](\cdot)  \, \, \mbox{on} \, \, \mathbb{R}.   $                                                                                                                                                                                                                          
 
 \item \emph{(Monotonicity in $M$)}\label{Thm:monotonicity_properties_1_to_0:prop_monotonicity_in_M}
Let  $0 \leq M <\widetilde{M}$ and $L \geq 0$ be fixed. Then $\mathscr{E}\left[\theta_{(-\widetilde{M},L)}^{(c)}\right](x) \leq  \mathscr{E}\left[\theta_{(-M,L)}^{(c)}\right](x).$

\item \emph{(Monotonicity in $L$)}\label{Thm:monotonicity_properties_1_to_0:prop_monotonicity_in_L}
Let  $0 \leq L <\widetilde{L}$ and $M \geq 0$ be fixed. Then $\mathscr{E}\left[\theta_{(-M,L)}^{(c)}\right](x) \leq  \mathscr{E}\left[\theta_{(-M,\widetilde{L})}^{(c)}\right](x).$
\item \emph{(Monotonicity in $x$)}\label{Thm:monotonicity_properties_1_to_0:1D_monotonicity_in_x} For every fixed $M$ and $L$ the mapping $x\mapsto \mathscr{E}\left[\theta_{(-M,L)}^{(c)}\right](x) $ is non-increasing.

\item \emph{(Continuous dependence of $\theta_{(-M,L)}^{(c)}(\cdot)$ on $L,M$)}\label{Thm:monotonicity_properties_1_to_0:1D_lemma_continuity_in_L_M} Let $0<M< \infty$, $0 <L< \infty$. The mappings $L\mapsto \mathscr{E}\left[ \theta_{(-M,L)}^{(c)}\right ](\cdot)$ and $M\mapsto \mathscr{E}\left[ \theta_{(-M,L)}^{(c)}\right ](\cdot)$ are continuous in the sup norm.
\end{enumerate}
\end{Lemma}

\subsection{Passing to the limit}
We are now ready to pass to the limit $M=\infty$ as a first step towards the proof of  Proposition \ref{DQ:c_zero:1_leadsto_0}. Define 
 \begin{eqnarray*}
\theta_{(-\infty,L)}^{(c)}(x) := \inf_{M >0}\mathscr{E}\left[\theta_{(-M,L)}^{(c)}\right](x)  = \lim_{M \to \infty}\mathscr{E}\left[\theta_{(-M,L)}^{(c)}\right](x),
 \end{eqnarray*}
where the last equality is a consequence of Lemma \ref{Thm:monotonicity_properties_1_to_0}(i). The following proposition highlights the role of the front speed $c$: roughly speaking it says that stretching procedure $M \to \infty$ we designed ``looses mass'' whenever $c>2$, i.e., $\theta_{(-\infty,L)}^{(c)}(x) \equiv 0$ when $c>2$.  Although zero is a (trivial) solution to the $\onelc$-truncated problem, one might wonder about the usefulness of the minimax construction we developed in \eqref{minimax}, for it seems to be not good  enough to obtain nontrivial solutions to \eqref{H_1D_problem:truncated} in $(-\infty,L)$. It turns out that the limitation is not on the method, but on the nature of the problem: no solution to \eqref{theorem:1D_result-equation} exists  when $c>2$, as we will show afterwards in Lem. \ref{1_leadsto_0_1D:existence}.  

\begin{Proposition}[Dichotomy $c \gtrless 2$]\label{Prop:Dichotomy} For any fixed $c>0$ we verify
\begin{align}\label{Prop:Dichotomy:equality}
 \theta_{(-\infty, L)}^{(c)}(\cdot) \neq 0, \quad \mbox{whenever} \quad c<2; & \qquad \theta_{(-\infty, L)}^{(c)}(\cdot) \equiv 0,\quad \mbox{whenever} \quad  \quad c>2.\end{align}
Furthermore, whenever $0 \leq c<2$ and $L>0$, the family $\theta_{(-\infty,L)}^{(c)}(\cdot)$ has the following properties: 
\noindent
 \begin{enumerate}[label=(\roman*), ref=\theTheorem(\roman*)]
 \item \label{Prop:Dichotomy:1D_continuity_in_L_M} The mapping $L\mapsto \theta_{(-\infty,L)}^{(c)}(\cdot)$ is continuous in the sup norm on $0\leq L\leq\infty$.
 \item {(Monotonicity)}\label{Prop:Dichotomy:1D_monotonicity_infty_L}
 The functions $x\mapsto \theta_{(-\infty,L)}^{(c)}(x)$ are defined for every $x \in \mathbb{R}.$
  The mapping  $L \mapsto \theta_{(-\infty,L)}^{(c)}(x)$ is non-decreasing for any fixed x. Furthermore, the mapping  $x \mapsto \theta_{(-\infty,L)}^{(c)}(x)$ is non-increasing  for any fixed L.
 \item\label{Prop:Dichotomy:minuns_infty_L_1D_theorem}The function $\theta_{(-\infty,L)}^{(c)}(x)$ solves \eqref{H_1D_problem:truncated} on  $(-\infty, L).$ Furthermore, $\displaystyle{\lim_{x \to -\infty}}\theta_{(-\infty,L)}^{(c)}(x) = 1.$
   \end{enumerate}   
 \end{Proposition}

\begin{Proof} We deal first with the case $c<2$: from ODE theory (cf. \cite[\S 4.4]{Fife}) there exists a solution $w(\cdot)\in \mathcal{C}^{\infty}(\R;[-1,1])$ to $\partial_x^2w + c\partial_xw + w - w^3=0$  satisfying $\displaystyle{\lim_{x\to-\infty}w(x) =1}$ and so that $w(x)$ is oscillatory as $x \to +\infty$ whenever $0 \leq c<2$; see Fig. \ref{Figure:1d-solutions_no_jump}. Translation invariance of solutions to this ODE allow us to assume without loss of generality that $ 0 = w(0) < w(x) <1$ for $x <0$. 

\begin{figure}[htb]
    \centering
    \begin{subfigure}[b]{0.47\textwidth}
          \centering
    \includegraphics[width=.9\textwidth,height=1.1in]{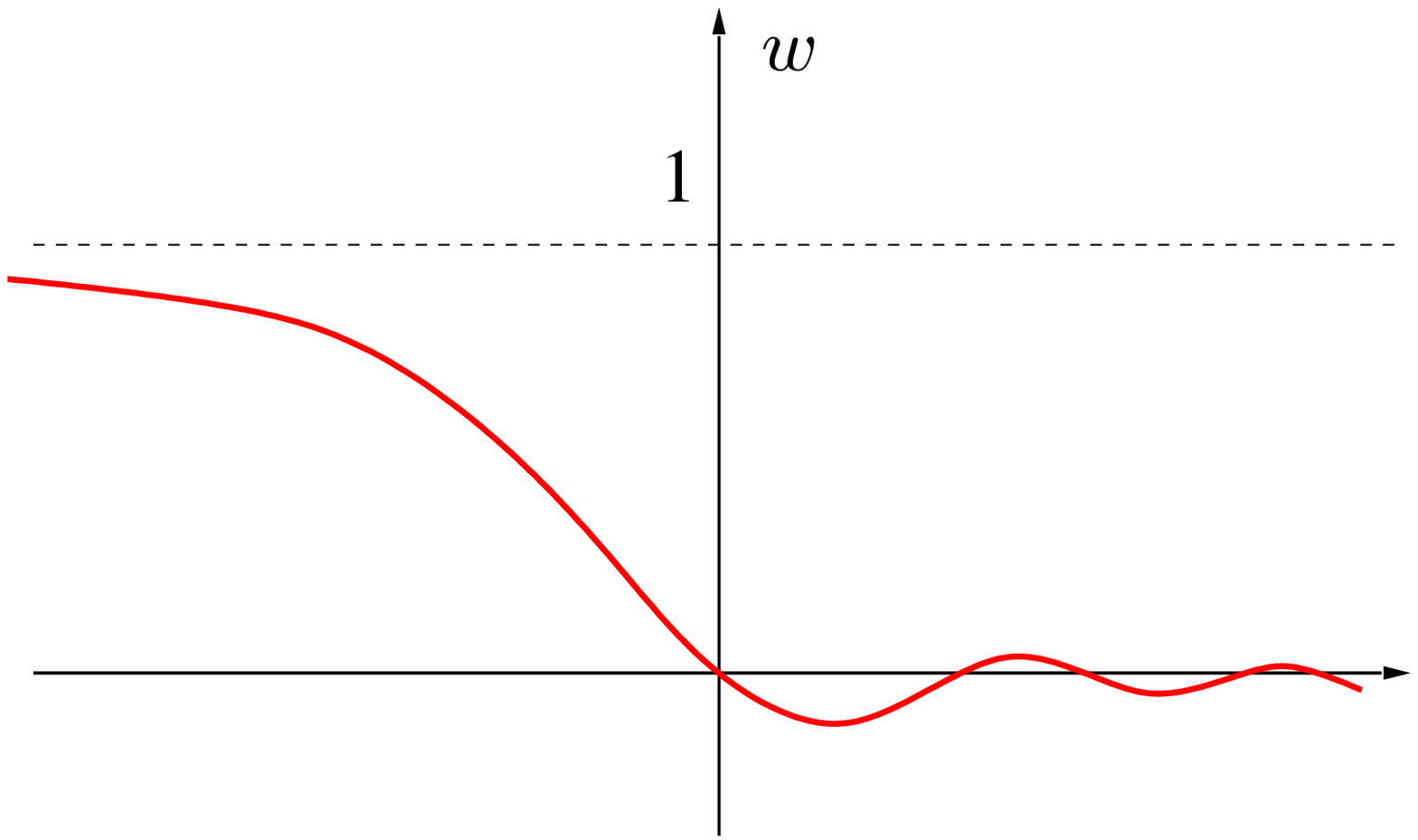}
    \end{subfigure}
    ~ 
    \begin{subfigure}[b]{0.47\textwidth}
    \centering    
          \includegraphics[width=.9\textwidth,height=1.1in]{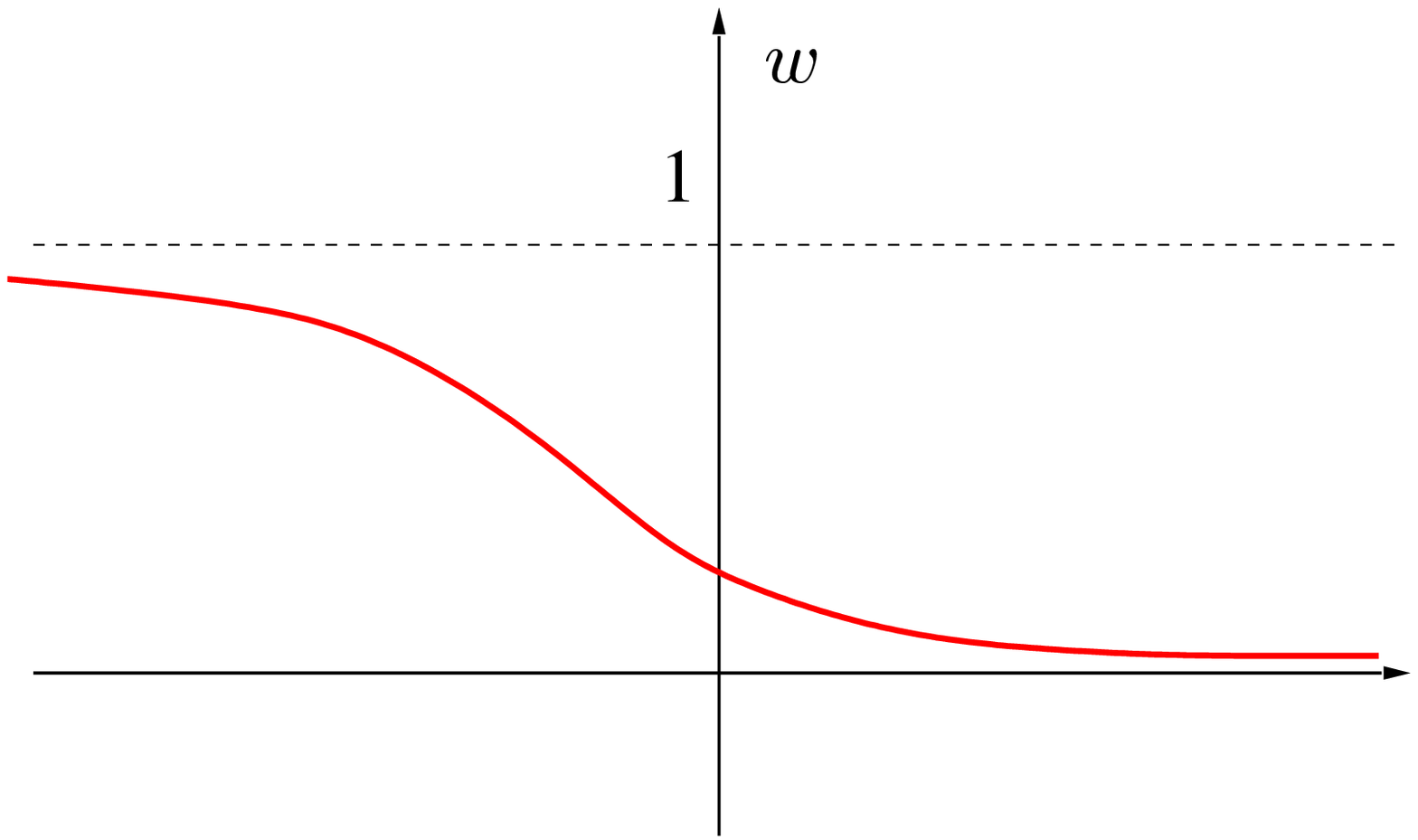}
        \end{subfigure}
    \caption{Sketch of solutions to $\partial_x^2w(x) + c\partial_xw(x) + w(x) - w^3(x) =0$ for $0\leq c<2$ (left) and $c>2$ (right) satisfying $\displaystyle{\lim_{x\to-\infty}w(x) = 1}$ and $\displaystyle{\lim_{x\to\infty}w(x) = 0}$.\label{Figure:1d-solutions_no_jump}}
\end{figure}
Applying classical comparison principles to the problem \eqref{H_1D_problem:truncated} on the interval $[-M,0]$ we conclude that $w(x) \leq \theta_{(-M,0)}^{(c)}(x)$ hence  $w(x) \leq \mathscr{E}\left[\theta_{(-M,0)}^{(c)}\right](x)\leq \mathscr{E}\left[\theta_{(-M,L)}^{(c)}\right](x)$ for $M>0$, thanks to Lem. \ref{Prop:Dichotomy:1D_continuity_in_L_M} and to the  monotonicity Lemma \ref{Prop:Dichotomy:1D_monotonicity_infty_L}. Taking the infimum in $M>0$ we conclude that $\theta_{(-\infty,L)}^{(c)}(\cdot)\neq 0$, which proves the first part of the result. As a byproduct we obtain (ii) using a squeezing property, for
$$\displaystyle{1 = \lim_{x \to -\infty}w(x) \leq \liminf_{x\to -\infty} \theta^{(c)}_{(-\infty,L)}(x)\leq 1}.$$ 
Item (i) is a direct consequence of Lem. \ref{Thm:monotonicity_properties_1_to_0:1D_lemma_continuity_in_L_M}. To show that $\displaystyle{\lim_{x\to\infty}\theta^{(c)}(x)=0}$ we use the function $\displaystyle{\bar{w}(x) := \frac{\csch(x + x_0)}{\sqrt{2}}}$  appropriately shifted so that $\bar{w}(0) =1$; notice that $\bar{w}$ satisfies $\partial_x^2\bar{w} - \bar{w} - (\bar{w})^3 =0$ and that is it monotonic, i.e., $\partial_xw(\cdot) \leq 0$. Hence, $w(\cdot)$ is a supersolution on any interval $[0,L]$ and classical comparison principles imply that 
$$ \theta^{(c)}_{(0,L)}(x) \leq \bar{w}(x), \quad \mbox{on} \quad x \in[0,L].$$
Thus,  $\theta_{(-\infty,L)}^{(c)}(x) \leq \bar{w}(x)$ using Lem. \ref{Thm:monotonicity_properties_1_to_0:prop_monotonicity_in_M}. The result is obtained using that $\theta^{(c)}(x) \geq 0$ and $\bar{w}(x) \to 0$ exponentially fast as $x \to \infty$. 
In order to prove the strict monotonicity of the solution, we use Prop. \ref{Thm:monotonicity_properties_1_to_0:1D_monotonicity_in_x}: the mapping $x\mapsto \theta^{(c)}(x)$ is monotonic in $x$ as the sup of monotonic functions, i.e., $\partial_x \theta^{(c)} \leq 0$. One obtains $\partial_x \theta^{(c)} < 0$ by applying Hopf lemma and the maximum principle (notice that the discontinuity of the control parameter $\mu(\cdot)$ plays no role here since, by classical regularity theory, we know that $\theta^{(c)}(\cdot)$ is in fact smooth away from the quenching front).

We now study the case $c>2$, showing that $\theta_{(-\infty,L)}^{(c)}(\cdot) \equiv 0$. We argue by  contradiction: assume the existence of a solution $\theta_{(-\infty,L)}^{(c)}(\cdot)$ satisfying \eqref{H_1D_problem:truncated} in $(-\infty,L)$ and so that $\theta_{(-\infty,L)}^{(c)}(x) >0$ when  $x \in (-\infty, L)$. There exists a $d \in (2, c)$ and a solution $v(\cdot)$ to $\partial_x^2v + d\partial_xv + v - v^3=0$  satisfying $v(-\infty) =1$, $v(\infty) =0$, $v(\cdot) >0$, $\partial_x v(\cdot) <0$. Define $m(x) = - \frac{x}{2}+ \sqrt{\frac{x^2}{4} +2}$. Results from the asymptotic theory of ODEs  (cf. \cite[Chap. 3, Sec. 8]{Coddington}) show that
\begin{align}\label{1_leadsto_0-truncated:asymptotics}
v(x) = 1 + \mathcal{O}(e^{m(d)x}),   \quad  \theta_{(-\infty,L)}^{(c)}(x) = 1 + \mathcal{O}(e^{m(c)x}), \quad \mbox{whenever} \quad x\to -\infty.
\end{align}
As $m(c) <m(d)$, there exists a $R>$ sufficiently large such that  $v(x) \geq  \theta_{(-\infty,L)}^{(c)}(x) $ for all $x \leq -R$. Positivity of $v(\cdot)$ implies  that  $v(x)\geq \theta_{(-\infty,L)}^{(c)}(x)$ whenever $x \geq L$. As $\theta_{(-\infty,L)}^{(c)}(\cdot)$ is monotone, there exists an $\epsilon >0$ such that 
$$\theta_{(-\infty,L)}^{(c)}(x) \leq 1 - \epsilon, \quad \mbox{for} \quad  x \in [-R,L].$$ In conclusion, we can make use of monotonicity of $v(\cdot)$ 
to obtain a $\tau \in \R$ such that  $w(x + \tau) \geq \theta_{(-\infty,L)}^{(c)}(x)$ on $x \in [-R,L]$, hence  guaranteeing that $v(x + \tau) \geq  \theta_{(-\infty,L)}^{(c)}(x)$ for all $x\in \R$. For such a $v(\cdot)$ we define $z(x) := v(x + \tau) - \theta_{(-\infty,L)}^{(c)}(x) \geq 0$. Using the properties of $v(\cdot)$ and the assumption  $\theta_{(-\infty,L)}^{(c)}(\cdot)\not \equiv 0$ we can find a shift $\tau$ in such a way that $z(\cdot)$ vanishes in at least one point. Properties of both $\theta_{(-\infty,L)}^{(c)}(\cdot)$ and $v(\cdot)$ imply that
\begin{equation}\label{w_is_supersolution}
 \partial_x^2 z(x) + c \partial_x z(x) + \mu(x)z(x) + f[\theta_{(-\infty,L)}^{(c)}(x),w(x)] z(x) = (c-d)\partial_x w(x) + [\mu(x) - 1]w(x)\leq 0,
\end{equation}
in $x\in (-\infty,L)$, where $f[a,b]:= \frac{a^3 - b^3}{a-b}$ whenever $a\neq b$ and $3a^2$ otherwise. One conclude that $z(\cdot)$ is a supersolution. As $z(\cdot) \geq 0$ has an interior minimum point, the maximum principle allied to the Hopf Lemma implies that $z \equiv 0$. However the latter is equivalent to $v(\cdot)\equiv\theta_{(-\infty,L)}^{(c)}(\cdot)$, which is a contradiction, for $v(\cdot)$ and $\theta_{(-\infty,L)}^{(c)}(\cdot)$ satisfy different equations. It finishes the proof. \end{Proof}


\begin{Lemma}[Existence/non-existence;   $\onelc$ problem]\label{1_leadsto_0_1D:existence} There exists a solution $\theta^{(c)}(\cdot) \in \mathcal{C}^{(1,\alpha)}(\R)$
 
\end{Lemma}

\begin{Proof}[of Theorem \ref{theorem:1D_result}] We begin by proving existence when $c<2$, following the ideas in \cite[Sec. 2]{Monteiro_Scheel} to which we refer to for further details: define $\displaystyle{\theta^{(c)}(x) = \sup_{L>0}\theta_{(-\infty,L)}^{(c)}(x) = \lim_{L\to \infty}\theta_{(-\infty,L)}^{(c)}(x)}.$ The asymptotic behavior at $x \to -\infty$ is a consequence of Prop. \ref{Prop:Dichotomy:1D_monotonicity_infty_L}-\ref{Prop:Dichotomy:minuns_infty_L_1D_theorem}, since  the mapping $L \mapsto \theta_{(-\infty,L)}^{(c)}(\cdot)$ is monotonic.

Next, we turn to the proof of  nonexistence when $c>2$: choosing $d$ so that $2<d<c$ we can find a strictly positive function $w(\cdot)$, $\partial_xw(\cdot) < 0$, satisfying 
\begin{equation}\label{w_properties}
 \partial_x^2w(x) + d\partial_xw(x) + w(x) - (w(x))^3 = 0, \qquad \lim_{x\to -\infty}w(x) = 1, \qquad \lim_{x\to -\infty}w(x) = 0.
\end{equation}
A direct computation shows that $\partial_x^2w(x) + c\partial_xw(x) + \mu(x)w(x) - (w(x))^3 = (\mu(x) -1)w +(c-d)\partial_x w\leq 0$. 
A similar analysis to that in \eqref{1_leadsto_0-truncated:asymptotics} shows that $w(x) \geq \theta^{(c)}(x)$ for all $x \to -\infty$. In order to understand and compare $\theta^{(c)}(x)$ and $w(x)$ as $x \to \infty$ we use an analysis similar to that of \cite{ber}: we have that  $w(\cdot),\, \theta^{(c)}(\cdot) >0$, and both  satisfy
\begin{align*}
 \partial_x^2w + d\partial_xw - w\leq 0 \leq \partial_x^2\theta^{(c)}+c\partial_x\theta^{(c)} - \theta^{(c)}.
\end{align*}
We conclude from \cite[Prop. 4.1 and Prop. 4.2]{ber} the existence of positive constants $M$ and $K$ such that 
\begin{align*}
 w(x)\geq K \exp\left(\frac{-d}{2} - \sqrt{ \frac{d^2}{4} -1}\right),\quad \theta^{(c)}(x)\leq M \exp\left(\frac{-c}{2} - \sqrt{ \frac{c^2}{4} +1}\right).
\end{align*}
Reasoning as in the proof of Prop. \ref{Prop:Dichotomy:equality} one obtain an $R>0$ so that 
\begin{align}\label{far_field_positivity}
w(x) - \theta^{(c)}(x) >0, \qquad \mbox{for} \quad x\in (-\infty,-R]\cup [R, \infty).
\end{align}
Now, as both $w(\cdot)$ and $\theta^{(c)}(\cdot)$ are bounded and non increasing, satisfying respectively the asymptotic properties \eqref{w_properties} and \eqref{H_1D_problem:truncated}, we conclude that we can shift $w(\cdot)$ so that $w(x - \tau) \geq \theta^{(c)}(x)$ for all $x \in \R$.
Now we take the infimum of $\tau$ such that $w(x - \tau) \geq \theta^{(c)}(x)$ holds with an equality in at least one point (clearly, the asymptotics in $w(\cdot)$ and $\theta^{(c)}(\cdot)$ shows that $\tau <+\infty$). Now, defining $z(\cdot) = w(\cdot - \tau) - \theta^{(c)}(\cdot)$ we get that z solves an inequality as that of \eqref{w_is_supersolution},
thanks to both the non increasing property of $w(\cdot)$ and its positivity. We conclude from the maximum principle and the Hopf's lemma that $z \equiv 0$, which is an absurd due to the asymptotic behavior of $w(\cdot)$ and $\theta^{(c)}(\cdot)$, and this contradiction gives the result.
\end{Proof}
From the properties of the subsolution used in the previous proof we readily derive the next result:
\begin{Lemma}[Exponential convergence; $\onelc$ problem]\label{1_leadsto_0_1D:exponential_convergence_ff} For any $\theta^{(c)}(\cdot)$ solving \eqref{theorem:1D_result-equation} when $c<2$ there exists a $C$, $\delta >0$ independent of $x$ such that
\begin{align*}
\lim_{x\to \infty}|\theta^{(c)}(x)| \leq Ce^{-\delta |x|}, \quad \mbox{and} \quad  \lim_{x\to \infty}|\theta^{(c)}(x) - 1 | \leq Ce^{-\delta |x|}.
\end{align*}
\end{Lemma}
To finalize this section we  show that the solutions obtained in \cite{Monteiro_Scheel} for small $c\geq0$ through perturbation methods agree with those constructed here. In passing we show that their continuity in the parameter $c$.
\begin{Lemma}[Uniqueness of the  continuation in $c$;   $\onelc$ fronts]\label{1_leadsto_0_1D:uniqueness} There exists a unique continuation $\theta^{(c)}(\cdot)$ solving the problem \eqref{1_leadsto_0:1D} whenever $c \in [0, 2)$. Moreover, the mappings $c\mapsto \theta^{(c)}(\cdot): [0,2) \to L^{\infty}(\R;\R)$ is continuous.
\end{Lemma}
\begin{Proof}
The proof consists in showing that for any $d\in[0,2)$ the linearization of the equation \eqref{Main_equation} at $\theta^{(d)}(\cdot)$ is invertible, i.e., 
\begin{align}\label{1_leadsto_0:linearization}
 \mathscr{L}_{\theta^{(d)}}[v] = \partial_x^2 v + d\partial_x v + \mu(x) v  - 3(\theta^{(d)})^2 v , \quad \mathcal{D}\left(\mathscr{L}_{\theta^{(d)}}\right) = H_2(\R).
\end{align}
is a boundedly invertible operator from $H^2$ to $L^2$. Indeed, assume the latter to be true. Then, plugging $\theta^{(d)} + u$ in \eqref{Main_equation} we rewrite it as
\begin{align*}
 \mathscr{L}_{\theta^{(d)}}[u] = \mathscr{N}[\theta^{(d)},u] + (d-c)\partial_x \theta^{(d)}.
\end{align*}
As $\mathscr{N}[\theta^{(d)},u] =\mathcal{O}(|u|^2)$ the term  on the right hand side is in $L^2(\R)$ and we can apply the IFT to obtain existence and  the uniqueness  of solutions in a neighborhood of $(\theta^{(d)}(\cdot), d)$. Furthermore, the mapping $d\mapsto \theta^{(d)}(\cdot)$ is continuous in $L^{\infty}(\R)$, thanks to the Sobolev embedding $H^2(\R) \hookrightarrow L^{\infty}(\R)$ (cf. \cite[\S 8 \& 9]{brezis}).
Keeping the previous discussion in mind, we devote the rest of the proof to showing the invertibility of operator $\mathscr{L}_{\theta^{(c)}}$. 

Initially, we write the conjugate operator  $\widetilde{\mathscr{L}_{\theta^{(c)}}}[\cdot] := e^{\frac{d\,x}{2}}\mathscr{L}_{\theta^{(c)}}[e^{\frac{-d\,x}{2}}\,\cdot]$:
\begin{align}\label{1_leadsto_0:linearization-self_adj}
 \widetilde{\mathscr{L}_{\theta^{(c)}}}[\tv] = \partial_x^2\tv + \left(\mu(x)-\frac{d^2}{4} \right)\tv  - 3(\theta^{(d)})^2 \tv , \quad \tilde{v} \in \mathcal{D}\left(\widetilde{\mathscr{L}_{\theta^{(d)}}}\right) = e^{\frac{d\,x}{2}}H^2(\R) =: H_{d}^2(\R).
\end{align}
 Defining $||u||_{H_d^2} = ||e^{\frac{d\,x}{2}}u||_{H^2}$, we see that the isometry $\mathcal{I}: H_d^2(\R) \to H^2(\R)$, $\mathcal{I}[u]= e^{\frac{d\,x}{2}}u$  implies that $\mathscr{L}_{\theta^{(d)}}$ is invertible on $H^2(\R)$ if and only if  $\widetilde{\mathscr{L}}_{\theta^{(d)}}$ is invertible on $H_c^2(\R)$. For a moment, consider the operator $\widetilde{\mathscr{L}_1}[\cdot] := \widetilde{\mathscr{L}_{\theta^{(d)}}}[\cdot]$ with domain  $\mathcal{D}\left(\widetilde{\mathscr{L}_1}\right)= H^2(\R)$; the analysis in \cite[\S 5]{Monteiro_Scheel} shows that this  is a self-adjoint, Fredholm operator of index $0$, with essential spectrum contained in $\{z\in \C| \mathrm{Re}(z) < 0\}$. In order to show invertibility we  show that this operator has a trivial kernel, which is proved as follows: the properties of the operator $\widetilde{\mathscr{L}_1}$ imply that the $\sigma\left(\widetilde{\mathscr{L}_1}\right)\cap \{x \in \R| x\geq 0\}$ is either empty or consists of point spectrum only. It is straightforward to show that this set is bounded, therefore assume that there exists a $\lambda_0 \geq 0$ maximal eigenvalue, with corresponding eigenfunction $u_0$. In the referred paper it was also proven that $u_0 \in Ker\left(\widetilde{\mathscr{L}_1}\right)$ is spatially localized, namely, 
 \begin{align}\label{spatial_localization:u_0}
|\nabla u_0(x,y)|+|u_0(x,y )| \leq  C e^{- \delta|x|} ,\quad  (x, y)  \in  \mathbb{R}\times [0,\kappa] \quad a.e.\qquad \mbox{whenever} \quad u_0 \in \mathrm{Ker}\left(\mathscr{L}_{\Xi}\right).  
 \end{align}
%
In fact, we know that we can take $\delta =\frac{d}{2}, $ thanks to the results in \cite[\S 4]{Monteiro_Scheel-contact_angle}. From the self-adjoint properties of $\widetilde{\mathscr{L}_1}$ we derive $u_0$ is a ground state associated to its maximal eigenvalue $\lambda_0\in \R$, therefore it satisfies $u_0(\cdot) \geq 0$  almost everywhere (cf. \cite[XII.12]{Reed_Simon}).  We can write the eigenvalue equation $\widetilde{\mathscr{L}_{1}}[u_0] = \lambda_0 u_0$ as
%
\begin{align}\label{eigenvalue_problem:u}
 \widetilde{\mathscr{L}_{1}}[u_0] = \partial_x^2 u_0(x) + \left(\mu(x)-\frac{d^2}{4} \right)u_0(x)  - 3e^{-d\,x}(\widetilde{\theta^{(d)}}(x))^2 u_0(x) = \lambda_0 u_0(x)
\end{align}
Setting  $\widetilde{\theta^{(d)}}(\cdot) =e^{\frac{d\, x}{2}} \theta^{(d)}(\cdot)$ and using the properties of the function $\theta^{(d)}(\cdot)$, we have
\begin{align}\label{1_leadsto_0_1D:rescaled_theta}
\partial_x^2 \widetilde{\theta^{(d)}}(x) + \left(\mu(x) - \frac{d^2}{4}\right) \widetilde{\theta^{(d)}}(x) - e^{-d\,x} (\widetilde{\theta^{(d)}})^3(x) = 0;
\end{align}
asymptotic theory of ODEs (cf. \cite[Chap. 3, Sec. 8]{Coddington}) implies that $\displaystyle{\lim_{|x|\to \infty}\widetilde{\theta^{(d)}}(x) =0}.$ Now, multiply \eqref{eigenvalue_problem:u} by $\widetilde{\theta^{(d)}}(\cdot)$ and \eqref{1_leadsto_0_1D:rescaled_theta} by $u_0(\cdot)$ subtract both equations and integrate in $\R$ to find
\begin{align*}
 \int_{\R} (\widetilde{\theta^{(d)}}(\mathrm{x}) \partial_x^2u_0(\mathrm{x}) - u_0(\mathrm{x})\partial_x^2\widetilde{\theta^{(d)}}(\mathrm{x}) d\mathrm{x} -2 \int_{\R} e^{-d\,x}\left(\widetilde{\theta^{(d)}}(\mathrm{x})\right)^2\widetilde{\theta^{(d)}}(\mathrm{x}) u_0(\mathrm{x}) d\mathrm{x} = \lambda_0 \int_{\R} \widetilde{\theta^{(d)}} u_0  d\mathrm{x}
\end{align*} 
Integration by parts shows that the first integral vanishes, thanks to the decay estimates for $\widetilde{\theta^{(d)}}(\cdot)$ and $u_0(\cdot)$. We are left with  
\begin{align*} 
-2 \int_{\R} e^{-d\,x}\left(\widetilde{\theta^{(d)}}(\mathrm{x})\right)^2\widetilde{\theta^{(d)}}(\mathrm{x}) u_0(\mathrm{x}) d\mathrm{x} = -2 \int_{\R} \left(\theta^{(d)}(\mathrm{x})\right)^2\widetilde{\theta^{(d)}}(\mathrm{x}) u_0(\mathrm{x}) d\mathrm{x}= \lambda_0 \int_{\R} w u_0  d\mathrm{x}.
\end{align*}
We observe that the spatial localization of $u_0(\cdot)$ as asserted in \eqref{spatial_localization:u_0} and the fact that $\widetilde{\theta^{(d)}}(x) = \mathcal{O}(e^{-\frac{d\,x}{2}})$ as $x \to -\infty$ imply that both integrals are finite. Sign considerations of both $\widetilde{\theta^{(d)}}(\cdot)$ and  $u_0(\cdot)$ show that the right-hand side is non-positive while the left-hand  is nonnegative (since $\lambda_0 \geq 0$), therefore the integral on the left is zero. Now, invoking the strict positivity of the pattern $\theta^{(d)}(\cdot)$ (or equivalently, that of $\widetilde{\theta^{(d)}}(\cdot)$) we conclude that $u_0(\cdot) \equiv 0$ almost everywhere, which contradicts the fact that $u_0(\cdot)$ is an eigenfunction. Therefore, no eigenvalue can be found on $\{z\in \mathbb{C}| \mathrm{Re}(z) \geq 0\}$; in other words, the operator $\widetilde{\mathscr{L}_{1}}$ is boundedly invertible. 

An intermediate step is necessary in order to go back to the operator $\widetilde{\mathscr{L}_{\theta^{(d)}}}$: first, define the family of weighted Sobolev spaces $H_{(\delta,\delta)}^2(\R) =  e^{\delta<x>}H^2(\R),$ where $<x> := \sqrt{1 + x^2}$. The action of the operator $\widetilde{\mathscr{L}_{1}}$ on these spaces can be studied  by the operators $\widetilde{\mathscr{L}_{1}^{(\delta)}}[\cdot] = e^{-\delta<x>}\widetilde{\mathscr{L}_{1}}[e^{\delta<x>}\cdot],$ defined as $H^2(\R) \to L^2(\R)$ mappings; we point out that the mapping $\delta \mapsto \widetilde{\mathscr{L}_{1}^{(\delta)}}$ is continuous in the operator norm. Standard Fourier analysis shows that the far field operators $\displaystyle{\widetilde{\mathscr{L}_{1}^{(\delta;\pm \infty)}} = \lim_{x\to \pm \infty}\widetilde{\mathscr{L}_{1}^{(\delta)}}}$ are boundedly invertible operators from  $H^2(\R)$ to $L^2(\R)$ for any $|\delta|\leq \frac{d}{2}$, hence  for any $\delta$ in this range the operators $\widetilde{\mathscr{L}_{1}^{(\delta)}}$ are Fredholm (cf. \cite[Prop. 4.3 and Rem. 4.7]{Monteiro_Scheel-contact_angle}). Continuity in $\delta$ implies that these operators also have index $0$. Therefore, proving invertibility is equivalent to showing that the kernel is trivial. In that regard, observe that we have a  scale of Banach spaces, i.e., $\displaystyle{H_{(\delta',\delta')}^2(\R) \subset H_{(\delta,\delta)}^2(\R)}$ whenever $\delta' > \delta$. Hence, one can invoke \cite[Lemma 5.3]{Monteiro_Scheel} or \cite[Lemma 4.6]{Monteiro_Scheel-contact_angle} to derive the persistence of elements in the kernel, namely, the equality $\displaystyle{\mathrm{Ker}\left(\widetilde{\mathscr{L}_{1}^{(\delta)}}\right) = \mathrm{Ker}\left(\widetilde{\mathscr{L}_{1}}\right)} =\{0\}$ holds for any $|\delta|\leq \frac{d}{2}$; being Fredholm operators of index 0 this property is equivalent to their invertibility.

We finally go back to the family of operators $\widetilde{\mathscr{L}_{\theta^{(d)}}}[\cdot]$. Clearly,  $\displaystyle{ e^{\frac{d}{2}x}H^2(\R)  \subset H_{(d,d)}^2(\R) = \mathcal{D}\left(Ker\left(\widetilde{\mathscr{L}_1}^{(\delta)}\right)\Big|_{\delta = \frac{d}{2}}\right)}$. Since the latter set is trivial, the same is also true of the kernel of $\widetilde{\mathscr{L}_1}$ taken with domain $e^{\frac{d}{2}x}H^2(\R)$, which corresponds to the operator $\widetilde{\mathscr{L}_{\theta^{(d)}}}[\cdot]$. As the latter is a Fredholm operator with index 0, this property immediately implies  bounded invertibility, and we conclude the proof.
\end{Proof}

\begin{Proof}[of Prop. \ref{theorem:1D_result}]
Combine the proofs of Lem. \ref{1_leadsto_0_1D:existence}, \ref{1_leadsto_0_1D:exponential_convergence_ff} and  \ref{1_leadsto_0_1D:uniqueness}.
\end{Proof}
\section{Two-dimensional quenched patterns -- periodic horizontal interfaces: \texorpdfstring{$\hp_{\kappa}$}{H-kappa} patterns, \texorpdfstring{$\pi <\kappa < \infty$}{pi<kappa<infty}}\label{Sec:2D_result:k_finite}
In this section we prove Theorem \ref{theorem:2D_result} in the case $\pi<\kappa<\infty$. As mentioned before, it is important in our approach that the nonlinearity is odd so that we can restrict the study of equation \eqref{Main_equation} to the stripe $(x,y) \in \R\times[0,\kappa]$; any solution $U(\cdot,\cdot)$ in $\R\times[0,\kappa]$  is extended to the whole plane $\R^2$ by successive reflections $U(x,-y)=-U(x,y)$ and $U(x,\kappa+y)=-U(\kappa-y)$. The method of proof  is similar to the one used in the $\onelc$ problem, although the  construction of subsolutions is more involved; we mostly follow the arguments in \cite[Sec. 3 \& 4]{Monteiro_Scheel} by truncating $\mathcal{H}_{\kappa}$-problem (with parameter $\pi <\kappa < \infty$) to a strip $ \mathcal{S}_{(-M,L)} := (-M,L)\times(0,\kappa)$, with  Dirichlet boundary conditions:
  \begin{align}\label{H_kappa:truncated}
\left\{ \begin{array}{ccc}
\Delta_{x,y} U +c\partial_x U+  \mu(x) U - U^3&=0, & (x,y)\in\mathcal{S}_{(-M,L)} ,\\
U-h_{(-M,L)}& =0, & (x,y)\in \partial \mathcal{S}_{(-M,L)},         
        \end{array}\right.
\end{align}
where $h_{(-M,L)}(x,y) := \theta_{(-M,L)}^{(c)}(x)\cdot\bu(y;\kappa)$, for $\theta_{(-M,L)}^{(c)}(\cdot)$ solution to the 1-D problem  \ref{H_1D_problem:truncated} and $\bu(\cdot, \kappa)$ given in \eqref{periodic_solutions}. We obtain a unique solution $\Xi_{(-M,L)}^{(c;\kappa)} \in \mathscr{C}^{(1,\alpha)}(\overline{\mathcal{S}_{(-M,L)}};[0,1]) $ by an  iteration scheme
\begin{eqnarray}\label{H_2D_stripe_existence_sequence}%
\left\{\begin{array}{c}
  -  \Delta \Xi_{n+1} - c\partial_x\Xi_{n+1} + 5 \Xi_{n+1} =  (5 +\mu(x)) \Xi_n - \Xi_n^3  \\
 (\Xi_{n+1} - h_{(-M,L)})\big|_{\partial \mathcal{S}_{(-M,L)}}=0
\end{array}\right.
\end{eqnarray}
where $\Xi_0(\cdot)$ is chosen in the class 
\begin{align*}
 \Psi_{\mathcal{H}_{\kappa}} &:= \left\{\Xi \in \mathscr{C}^{(1,\alpha)}(\mathcal{S}_{(-M,L)})\,\middle|\,(\Xi_{0} - h_{(-M,L)})\big|_{\partial \mathcal{S}_{(-M,L)}} \geq 0,\ \Xi(x)\in[0,1] \mbox{ for all } x,y \right\}\\
 & \quad \cap \{  \Delta\Xi + \mu(x)\Xi- \Xi^3 \leq 0, \, \  \Delta\Xi + \mu(x)\Xi - \Xi^3 \not\equiv 0 \, \mbox{in the sense of distributions}\}.
\end{align*}
Throughout this section, we fix $\kappa \in (\pi, \infty)$ and suppress the dependence of $\Xi$ and $\bar{u}$ on $\kappa$ for ease of notation. As in the previous section, proofs that are similar to those in \cite{Monteiro_Scheel} are only outlined and details are referred to that paper.

\begin{Proposition}[Existence and uniqueness; truncated $\hp_{\kappa}$ problem]\label{H_kappa:existence_uniqueness} Problem \eqref{H_kappa:truncated} has a unique solution 
$\displaystyle{\Xi_{(-M,L)}^{(c)}(\cdot,\cdot) \in \mathscr{C}^{(1,\alpha)}([-M,L]\times [0,\kappa]; [0,1])}$, $\forall \alpha \in [0,1)$.\end{Proposition}
\begin{Proof} The existence is obtained as in  \cite[Lem. 3.1]{Monteiro_Scheel} using the iterative scheme \eqref{H_2D_stripe_existence_sequence}; uniqueness follows as in \cite[Prop. 3.2]{Monteiro_Scheel} and integration by parts, as in Lem. \eqref{Existence_uniqueness_truncated:1_leadsto_0}. The stated regularity is derived using classical results in elliptic theory, as shown in \cite[\S 3]{Monteiro_Scheel}.\end{Proof}

 We define extension operators in order to compare solutions for different values of $M,L$, namely,  
\begin{equation}\label{extension_operator:H_kappa}
 \mathscr{E}\left[\Xi_{(-M,L)}^{(c)}\right](x,y) = \left \{\begin{array}{cc}
                          \Xi_{(-M,L)}^{(c)}(x,y), &   \quad  \mbox{for} \quad  (x,y) \in \mathcal{S}_{(-M,L)}\\
                           \bu(y) ,&\quad  \mbox{for} \quad (x,y) \in (-\infty, -M)\times [0,\kappa]\\
                            0 ,&\quad  \mbox{for} \quad (x,y) \in (L, \infty)\times [0,\kappa],
                          \end{array}\right.   \end{equation} 
where $\bu(\cdot) = \bu(\cdot;\kappa)$ is given in \eqref{periodic_solutions}. We use the same symbols for the one- and two-dimensional extension operators, slightly abusing notation, distinguishing between the two through the domain of definition of the function $\mathscr{E}$ is applied to. 
The proofs of the following Lem. \ref{Thm:sub_super_solutions-H_kappa} and Prop. \ref{Thm:monotonicity_properties_H_kappa_problem} are obtained as  in \cite[\S 3]{Monteiro_Scheel}: 

\begin{Lemma}[Comparison principles; $\mathcal{H}_{\kappa}$-problem]\label{Thm:sub_super_solutions-H_kappa} Let $\Xi_{(-M,L)}^{(c)}(\cdot,\cdot)$ be the solution from Proposition \ref{H_kappa:existence_uniqueness}.
\noindent
\begin{enumerate}[label=(\roman*), ref=\theTheorem(\roman*)]
\hangindent\leftmargin
 \item \emph{(2D supersolutions)}\label{Thm:sub_super_solutions-H_kappa:2D_supersolutions} 
 If $v$ satisfies, in the sense of distributions,
\begin{equation*}
\Delta v  + c\partial_x v +  \mu(x) v- v^3 \leq 0,\ (x,y)\in \mathcal{S}_{(-M,L)}, \quad \left(v - h_{(-M,L)}\right)\big|_{\partial \mathcal{S}_{(-M,L})} \geq 0, \quad   0 \leq v \leq 1,
\end{equation*}
then $v \geq \Xi_{(-M,L)}^{(c)}$ in $\mathcal{S}_{(-M,L)}.$ In particular, $v \geq \Xi_{(-M,L)}^{(c)}$ in $\mathcal{S}_{(-M,L)}$ for any solution of 
\begin{equation*}
\Delta v  + c\partial_x v + \mu(x) v- v^3 =0,\ (x,y)\in \mathcal{S}_{(-M,L)}, \quad \left(v - h_{(-M,L)}\right)\big|_{\partial \mathcal{S}_{(-M,L})} \geq 0, \quad   0 \leq v \leq 1.
\end{equation*}

\item\emph{(2D subsolutions)}\label{Thm:sub_super_solutions-H_kappa:2D_subsolutions} 
If $v$ satisfies, in the sense of distributions,
\begin{equation}\label{Thm:sub_super_solutions-H_kappa:2D_supersolutions:supersolution_lemma_number_2-H_kappa}
\Delta v  + c\partial_x v + \mu(x) v- v^3 \geq 0,\ (x,y)\in \mathcal{S}_{(-M,L)}, \quad \left(v - h_{(-M,L)}\right)\big|_{\partial \mathcal{S}_{(-M,L})} \geq 0, \quad   0 \leq v \leq 1,
\end{equation}
then $v \leq \Xi_{(-M,L)}^{(c)}$ in $\mathcal{S}_{(-M,L)}.$ In particular, $v \leq \Xi_{(-M,L)}^{(c)}$ in $\mathcal{S}_{(-M,L)}$ for any solution of 
\begin{equation*}
\Delta v  + c\partial_x v +  \mu(x) v- v^3 =0,\ (x,y)\in \mathcal{S}_{(-M,L)}, \quad \left(v - h_{(-M,L)}\right)\big|_{\partial \mathcal{S}_{(-M,L})} \leq 0, \quad   0 \leq v \leq 1.
\end{equation*}

\item \label{Thm:sub_super_solutions-H_kappa:2D_supersolutions_M_L}
The functions $\bar{u}(\cdot;\kappa) $ and $\theta_{(-M,L)}^{(c)}(x)$ given by \eqref{periodic_solutions} and Theorem \ref{theorem:1D_result}, respectively, are supersolutions to \eqref{H_kappa:truncated} in $\mathcal{S}_{(-M,L)}$. Consequently, for any $L, M >0$ we have that 
 \begin{eqnarray*}
 \Xi_{(-M,L)}^{(c)}(x,y)\leq \min\{\theta_{(-M,L)}^{(c)}(x), \bar{u}(y)\}. 
 \end{eqnarray*}
 \end{enumerate}
\end{Lemma}

\begin{Proposition}[Properties of the extension operator, $\mathcal{H}_{\kappa}$-problem]\label{Thm:monotonicity_properties_H_kappa_problem} The following properties of  $\mathscr{E}[\cdot]$ hold.
\begin{enumerate}[label=(\roman*), ref=\theTheorem(\roman*)]
\hangindent\leftmargin
 \item \emph{(Monotonicity of $\mathscr{E}$)}\label{Thm:monotonicity_properties_H_kappa_problem:monotonicity_extension_operator}
  We have $0 \leq \mathscr{E}\left[\Xi_{(-M,L)}^{(c)}\right](\cdot, \cdot) \leq 1$. Furthermore, if $w$ is only defined in a subset $A \subset \mathbb{R}^2$ so that  $ \mathcal{S}_{(-M,L)} \subset A \subset \Omega_{(-\infty,L)}$, $0 \leq w(\cdot) \leq 1$ and $w(\cdot) \leq \Xi_{(-M,L)}^{(c)}(\cdot, \cdot)$,  then 
$0 \leq \mathscr{E}\left[w\right](\cdot) \leq \mathscr{E}\left[\Xi_{(-M,L)}^{(c)}(\cdot, \cdot)\right](\cdot)  \, \, \mbox{in} \, \, \mathbb{R}^2 $.                                                                                                                                                                                                                                        
 \item \emph{(Monotonicity in $M$)}\label{Thm:monotonicity_properties_H_kappa_problem:2D_monotonicity_in_M}
Let  $0 \leq M <\widetilde{M}$ and $L \geq 0$ be fixed. Then $
M <\widetilde{M} \implies   \mathscr{E}\left[\Xi_{(-\widetilde{M},L)}\right](x,y) \leq  \mathscr{E}\left[\Xi_{(-M,L)}^{(c)}\right](x,y). $
 \item \emph{(Monotonicity in $L$)} \label{Thm:monotonicity_properties_H_kappa_problem:2D_monotonicity_in_L}
Let  $0 \leq L <\widetilde{L}$ and $M \geq 0$ be fixed. Then $
L <  \widetilde{L} \implies   \mathscr{E}\left[\Xi_{(-M,L)}^{(c)}\right](x,y) \leq  \mathscr{E}\left[\Xi_{(-M,\widetilde{L})}\right](x,y). $
 \item \emph{(Monotonicity in $x$)}\label{Thm:monotonicity_properties_H_kappa_problem:2D_monotonicity_in_x}
Let  $L,M,y$be fixed. Then
the mapping  $ x \to \mathscr{E}\left[\Xi_{(-M,L)}^{(c)}\right](x,y)$ is non-increasing. 
\end{enumerate}
\end{Proposition}
\subsection{Passing to the limit}\label{properties_theta_minusinfty_L}

We are now ready to pass to the limit $M=\infty$. Define 
\begin{equation}\label{definition:theta_infty_L:H_kappa}
\Xi_{(-\infty,L)}^{(c)}(x,y) := \inf_{M >0 } \mathscr{E}\left[\Xi_{(-M,L)}^{(c)}\right](x,y) =\lim_{M \to +\infty } \mathscr{E}\left[\Xi_{(-M,L)}^{(c)}\right](x,y), 
\end{equation}
where the last equality holds due to monotonicity of the mapping $M \mapsto \Xi_{(-M,L)}^{(c)}(x,y)$, Proposition \ref{Thm:monotonicity_properties_H_kappa_problem}(ii).

In this section we verify monotonicity properties and limits at spatial infinity of the limits $\Xi_{(-\infty,L)}(x,y)$ constructed in \eqref{definition:theta_infty_L:H_kappa}.

 \begin{Lemma}[Monotonicity of $\Xi_{(-\infty,L)}(x,y)$]\label{H_kappa_L_problem_solution}
 The following properties hold:
 \begin{enumerate}[label=(\roman*), ref=\theTheorem(\roman*)]  
\hangindent\leftmargin
 \item \label{H_kappa_L_problem_solution:solution_in_hald_truncated_domain}The function $\Xi_{(-\infty,L)}^{(c)}(\cdot,\cdot)$ solves the problem \eqref{H_kappa:truncated} in  $\mathcal{S}_{(-\infty,L)}:= \{(x,y) \in \mathbb{R}^2 | x <L, y \in [0,\kappa] \}$. 
 \item  The function $\Xi_{(-\infty,L)}^{(c)}(\cdot,\cdot)$  is non-increasing in $x$ and non-decreasing in $L$. 
 \item \label{H_kappa_L_problem_solution:supersolution} The inequality $ \Xi_{(-\infty,L)}(x,y)\leq \min\{\theta_{(-\infty,L)}^{(c)}(x), \bar{u}(y)\} $
 holds for all $(x,y) \in \mathcal{S}_{(-M,L)},$ where  $\theta_{(-\infty,L)}^{(c)}(\cdot)$ is given by \eqref{theorem:1D_result} and $\bar{u}(\cdot)$ by Prop. \ref{periodic_solutions}.  In particular, we have  $\displaystyle{ \sup_{L>0}\left(\Xi_{(-\infty,L)}(x,y) \right)\leq \min\{\theta^{(c)}(x), \theta_{(-\infty,0)}^{(c)}(y)\}}.$
 \end{enumerate}

 \end{Lemma}

\begin{Proof}Assertion (i) and (ii) follow as in \cite[Lem. 3.6]{Monteiro_Scheel}, being a consequence of a comparison method and monotonicity  of $\Xi_{(-M,L)}^{(c)}(\cdot)$ and $\Xi_{(-\infty,L)}^{(c)}(\cdot)$ in their arguments. We derive the inequalities in (iii) by passing to the limit $L=\infty$ and using the fact that for all $x\in \R$ the mapping $L \mapsto \theta_{(-\infty,L)}^{(c)}(x)$ is non-decreasing. 
 \end{Proof}

 \begin{Remark}
  One can readily conclude from Lem. \ref{H_kappa_L_problem_solution:supersolution} that $\Xi_{(-\infty,L)}^{(c)} \equiv 0$ whenever $c >2$, $L>0$. In fact, the non-existence of nontrivial solutions happens in a wider range for the  parameter $c$, as we show next.
 \end{Remark}


 \subsection{Existence for the \texorpdfstring{$\mathcal{H}_{\kappa}$}{H}-problem: case \texorpdfstring{$\frac{c^2}{4} +\frac{\pi^2}{\kappa^2}<1$}{c2/4 + pi2/kappa2 <1} }

 In order to prove the existence of solutions we construct appropriate subsolutions with the help of the next lemmas:
\begin{Lemma}[Properties of the family of periodic solutions $\bar{u}(\cdot,\kappa)$]\label{Lemma:periodic-solutions} Let $\bar{u}(\cdot;\kappa)$ be a solution to \eqref{periodic_solutions} and $\kappa > \pi$.  The following two properties hold:
\begin{enumerate}[label=(\roman*), ref=\theTheorem(\roman*)]
\hangindent\leftmargin
 \item\label{Lemma:periodic-solutions:item_a}	The quantity $\displaystyle{M := \sup_{y\in [0,\kappa]}\bar{u}(y;\kappa) }$ satisfies $\displaystyle{ M < \sqrt{\left(1 - \frac{\pi^2}{\kappa^2}\right)}};$
 \item\label{Lemma:periodic-solutions:item_b} For any $0\leq \alpha \leq  M$ we have  $v(y):= \alpha \sin\left(\frac{\pi\,y}{\kappa}\right) \leq \bar{u}(y;\kappa),$ for all $y \in [0,\kappa].$
\end{enumerate}

\end{Lemma}
\begin{Proof}
To prove the estimate in (i) we use the elliptic integral that gives the relation between amplitude and spatial period given in \cite{Monteiro_Scheel}[Lemma 4.1, equation (4.4)], \cite[\S V]{Hale}, 
\begin{equation}\label{period_of_periodic_orbits}
 \kappa := \mathscr{\kappa}(M) = 2 \sqrt{2} \int_{0}^1 \frac{dv}{\sqrt{[(1 - v^2)(2 - M^2(1 + v^2))]}} = \frac{2 \sqrt{2} \gamma}{M} \int_{0}^1 \frac{dv}{\sqrt{[(1 - v^2)(1 - (\gamma\,v)^2)]}},
\end{equation}
for $\gamma^2 = \frac{M^2}{2 - M^2}$. Notice that $0 \leq \gamma<1$. We find a lower bound to the integral on the right hand side:
\begin{align*}
\mathscr{\kappa}(M) >  \frac{2 \sqrt{2} \gamma}{M} \int_{0}^1 \frac{dv}{\sqrt{[(1 - v^2)(1 - \gamma^2)]}} =  \frac{2 \sqrt{2} \gamma}{M\sqrt{1 - \gamma^2}} \int_{0}^1 \frac{dv}{\sqrt{(1 - v^2)}}=  \frac{\pi \sqrt{2} \gamma}{M\sqrt{1 - \gamma^2}}.
\end{align*}
Squaring both sides and plugging $\gamma$ we obtain $\displaystyle{\kappa^2 > \frac{\pi^2}{1 - M^2} \Longleftrightarrow M^2 < 1 - \frac{\pi^2}{\kappa^2}}$, which finishes the proof of (i). In order to prove (ii) we  exploits the structure of this ODE in \eqref{periodic_solutions}, whose Hamiltonian is $\displaystyle{\mathcal{H}(u, u_y) = u_y^2 + u^2 - \frac{u^4}{2}}$,  cf. \cite[\S 2.1]{Monteiro_Scheel}. Indeed, considering $\bar{u}(\cdot)$ a periodic orbit with period $\kappa$ and maximal amplitude $M$, we readily obtain that $\mathcal{H}(\bar{u}, \bar{u}_y) =   \frac{M^2(2 - M^2)}{2}$.
Let  $v(y) := \alpha \sin\left(\frac{\pi\, y}{\kappa}\right)$ and $z(x) := \bar{u}(y) - v(y)$. Whenever $0 \leq \alpha\leq M$ one can see that 
\begin{align}\label{compare_hamiltonians}
\mathcal{H}(v, v_y)\leq  \mathcal{H}(\bar{u}, \bar{u}_y).  
\end{align}
As $\bar{u}_y(0) >v_y(0)$ and $\bu(0) = v(0)=0$ it is clear that $z(x) >0$ for and $x>0$ sufficiently small. By translation invariance of the solutions to the ODE \eqref{periodic_solutions},  reversibility of the solutions $\bar{u}$ with respect to $x \mapsto -x$, and the fact that the mapping $y \mapsto  \sin(y + \pi/2)$ is even it suffices to show that $z(y) \geq 0$ for $0 \leq y \leq \frac{\kappa}{2}.$ Assume that there exists a $0<x_0<\frac{\kappa}{2}$ such that $z(x_0) =0$. As $z(\cdot)\geq 0$ solves the elliptic differential equation, we can find and $A>0$ sufficiently large so that 
$$\partial_y^2z(y) + (1 - f[\bar{u},v](y) - A)z(y) \leq \partial_y^2z(y) + (1 - f[\bar{u},v](y) )z(y)=0, $$
where $f[a,b]:= \frac{a^3 - b^3}{a-b}$ whenever $a\neq b$ and $3a^2$ otherwise. We conclude using Hopf's lemma that $\partial_y z(x_0) <0$, which is absurd, since the inequality \eqref{compare_hamiltonians} prevents it from happening. Therefore $z(y) \geq 0$, hence $v(y) \leq \bar{u}(y)$ for  $y \in \left[0,\frac{\kappa}{2}\right]$, and by symmetry, for $y \in \left[0,\kappa\right].$ \end{Proof}
\begin{Lemma}[Sub and supersolutions]\label{Lemma:H_kappa:sub} Choose $d \in (c,2)$  so that $\displaystyle{\frac{c^2}{4} + \frac{\pi^2}{\kappa^2} < \frac{d^2}{4}< 1}$ and let 
$ w(\cdot)\in \mathcal{C}^{\infty}(\R;[-1,1])$ be a solution to  $\partial_x^2w + d\partial_xw + w - w^3=0$  satisfying $w(-\infty) =1$ and so that $ 0 = w(0) < w(x) <1$ for $x <0$ (see. Prop. \ref{Prop:Dichotomy} and Fig. \ref{Figure:1d-solutions_no_jump}). Let $M>0$, $\alpha$ and $\bar{u}(\cdot)$ be given as in Lem. \ref{Lemma:periodic-solutions}. Define 
$$  V(x,y) := e^{-\frac{(c-d)\, x}{2}}w(x)\cdot v_{\kappa}(y), \qquad \mbox{where} \quad v_{\kappa}(y) := \alpha \sin\left(\frac{\pi\, y }{\kappa}\right).$$
 Then, 
\begin{equation}\label{subsolution_inequality}
V(x,y) \leq \Xi_{(-\infty,L)}^{(c)}(x,y) \quad \mbox{in} \quad  [-M,0]\times[0,\kappa].
\end{equation}
Furthermore, whenever  $\theta^{(c)}(\cdot)$ is given by \eqref{theorem:1D_result} and $\bar{u}(\cdot)$ by Prop. \ref{periodic_solutions}, the following inequality holds
\begin{eqnarray}\label{Thm:stretched:2D_stripe_supersolutions_infty_L-2}
 \Xi_{(-\infty,L)}^{(c)}(x,y)\leq \min\{\theta^{(c)}(x), \bar{u}(y)\}. 
 \end{eqnarray}  
\end{Lemma}
\begin{Proof} Inequality \eqref{Thm:stretched:2D_stripe_supersolutions_infty_L-2} is readily derived  from Lem. \ref{Thm:sub_super_solutions-H_kappa:2D_supersolutions_M_L}, taking infimum in $M>0$, using the definition of $\Xi_{(-\infty,L)}^{(c)}(\cdot,\cdot)$ (see \eqref{definition:theta_infty_L:H_kappa}) and  the monotonicity of the mapping  $L \mapsto \theta_{(-\infty,L)}^{(c)}(\cdot)$.

In order to obtain \eqref{subsolution_inequality} we use Lem. \ref{Thm:sub_super_solutions-H_kappa:2D_subsolutions}, showing that $V(\cdot,\cdot)\leq \Xi_{(-M,0)}^{(c)}(\cdot,\cdot)$ in  $\mathcal{S}_{(-M,0)} = [-M,0]\times[0,\kappa]$ for any fixed $M>0$. Thanks to Lemma \ref{Lemma:periodic-solutions:item_b} it is straightforward to show that   $V(x,y) \leq \Xi_{(-M,0)}^{(c)}(x,y)$ on the boundary $\partial \mathcal{S}_{(-M,0)}$. Now we show that V satisfies
\begin{align*}
 \Delta_{x,y}V + c\partial_x V + V - V^3 \geq 0,
\end{align*}
in $\mathcal{S}_{(-M,0)} = [-M,0]\times[0,\kappa]$. Indeed, a direct calculation shows that 
 \begin{align*}
\Delta_{x,y}V + c\partial_x V + V - V^3 &= V\left[\frac{d^2}{4} - \frac{c^2}{4}  - \frac{\pi^2}{\kappa^2}\right] + e^{-\frac{(c-d)\,x}{2}}w^3(x)\left[v_{\kappa}(y) - e^{-(c-d)\,x}v_{\kappa}^3(y)\right]  \geq 0,
 \end{align*}
since $V(\cdot,\cdot)$ is non-negative, $x\leq 0$ and $d >c$. It follows that $V(x,y) \leq \Xi_{(-M,0)}^{(c)}(x,y).$  We obtain \eqref{subsolution_inequality} by invoking the  monotonicity of the mapping $L \mapsto \Xi_{(-M,L)}^{(c)}(\cdot,\cdot)$ and using the definition \eqref{definition:theta_infty_L:H_kappa}. 
\end{Proof}

\begin{Lemma}[Existence; $\mathcal{H}_{\kappa}$-problem, $\pi <\kappa <\infty$]\label{H_kappa:existence}
 Equation $\eqref{Main_equation}$ has  a solution $\Xi^{(c)}(\cdot,\cdot)\in \mathscr{C}^{(1,\alpha)}\left(\R\times [0,\kappa];\R\right)$, for any $0\leq \alpha<1$, where the latter is defined as
 \begin{align}\label{H_kappa:solution}
 \Xi^{(c)}(x,y) := \lim_{L\to \infty} \Xi_{(-\infty,L)}^{(c)}(x,y)
\end{align}

\end{Lemma}

\begin{Proof}
Most of the proof goes as in the paper \cite[Prop. 3.11]{Monteiro_Scheel}. The monotonicity properties of the functions $\Xi^{(c)}_{(-\infty,L)}(\cdot,\cdot)$ show that the definition \eqref{H_kappa:solution} makes sense and Lem. \ref{Lemma:H_kappa:sub} shows that $\Xi_{(-\infty,L)}^{(c)}(\cdot,\cdot)\not \equiv 0$ whenever $\frac{c^2}{4} +\frac{\pi^2}{\kappa^2}<1$. As $0\leq \Xi_{(-\infty,L)}^{(c)}(\cdot,\cdot) \leq \Xi_{\kappa}^{(c)}(\cdot,\cdot)$ we conclude that $\Xi_{\kappa}^{(c)}(\cdot,\cdot)$ is also nontrivial. Using Lebesgue dominated convergence we conclude that 
exists in the pointwise sense and that this sequence converges in $L_{\mathrm{loc}}^1$ hence in the sense of distribution, solving the equation \eqref{Main_equation} in the domain $\R\times [0,\kappa]$. Now it remains to show that the asymptotic limits are satisfied, namely, that 
\begin{equation*}
 \lim_{x \to -\infty}\Xi_{\kappa}^{(c)}(x,y) = \bar{u}(y), \qquad  \lim_{x \to +\infty}\Xi_{\kappa}^{(c)}(x,y) = 0.
\end{equation*}
The limit on the right follows easily from inequality \eqref{Thm:stretched:2D_stripe_supersolutions_infty_L-2}, for $\displaystyle{\lim_{x\to \infty}\theta^{(c)}(x) =0}.$
The proof of the limit on the left is more involved, and our analysis has some similarities to those of \cite{Vega} and \cite[Theorem 1]{Wang_Hui}. %
Indeed, monotonicity results derived in Lem. \ref{Thm:monotonicity_properties_H_kappa_problem} allow us to conclude that 
\begin{equation}\label{far_field_x_minus_infty_ineq}
 v_L(y) := \lim_{x\to\infty }\Xi_{(-\infty,L)}^{(c)}(x,y)\leq \liminf_{x\to\infty }\Xi_{\kappa}^{(c)}(x,y)\leq \limsup_{x\to\infty }\Xi_{\kappa}^{(c)}(x,y)\leq \bar{u}(y),
\end{equation}
where the first limit is known to exists thanks to the monotonicity in $x$ of $\Xi_{(-\infty,L)}^{(c)}(\cdot, \cdot)$. According to Lem. \ref{H_kappa_L_problem_solution:solution_in_hald_truncated_domain}  we know that $\Xi_{(-\infty,L)}^{(c)}(\cdot, \cdot)$ solves the equation \eqref{Main_equation}. Thus, $v_L(\cdot)$   satisfies 
\begin{equation*}\label{limit_in_x-equation:L}
\partial_y^2v_L(y) + v_L(y) - \left(v_L(y)\right)^3=0
\end{equation*}
in the sense of distributions in $[0,\kappa]$, hence in the classical sense.  As $v_L(y)\Big|_{y = 0, \kappa} =0$ we conclude that either $v\equiv0$ or $v$ is a periodic solution with period $\tau$ so that $\displaystyle{\frac{\kappa}{\tau} \in \mathbb{N}}.$
We can readily exclude the first possibility, since Lem. \ref{Lemma:periodic-solutions:item_b} implies that  $\displaystyle{\alpha \sin\left(\frac{\pi\, y}{\kappa}\right)\leq v_L(y)\leq \bu(y)}$. The same inequality also implies that $\tau = 2\kappa$, i.e., $\bar{u}(\cdot)$ and $v_L(\cdot)$ have the same period therefore and obey the same normalization, therefore $v_L(\cdot) \equiv \bar{u}(\cdot)$, and the result follows from the equality of  \eqref{far_field_x_minus_infty_ineq}.\end{Proof}

Unlike in the previous case, it is not directly clear that the convergence $\displaystyle{\lim_{x\to -\infty}\Xi_{\kappa}^{(c)}(x,y) = \bar{u}(y)}$ has exponential rate of convergence. Our next result implies that.

\begin{Lemma}[Exponential convergence; $\mathcal{H}_{\kappa}$-problem, $\pi <\kappa <\infty$]\label{H_kappa:exponential_convergence} There exists a $C$, $\delta >0$ independent of $x$ and $y$ such that
\begin{align}
\lim_{x\to \infty}|\Xi_{\infty}^{(c)}(x,y)| \leq Ce^{-\delta |x|},\mbox{and} \lim_{x\to -\infty}|\Xi_{\infty}^{(c)}(x,y) - \bar{u}(x)|\leq Ce^{-\delta |x|}.
\end{align}
\end{Lemma}
\begin{Proof} 
Initially we show exponential rate of convergence to the far-field as $|x|\to \infty$. The result follows if we show that $\partial_x(\Xi_{\infty}^{(c)})\in e^{-\delta<x>}H^2(\R\times[0,\kappa])$ for some $\delta >0$ and $<x> :=\sqrt{1 +x^2}$. Indeed, as we know from Lem. \ref{H_kappa:existence}, $\displaystyle{\lim_{x\to -\infty}\Xi_{\infty}^{(c)}}(x,y) = \bar{u}(y)$; using the Sobolev embedding $H^{2}(\R\times[0,\kappa])\hookrightarrow  L^{\infty}(\R\times[0,\kappa])$ we have
\begin{align*}
 |\Xi_{\infty}^{(c)}(x,y) - \bar{u}(y)| \leq \left|\int_{\infty}^x \partial_x(\Xi_{\infty}^{(c)})(s,y)ds\right| \lesssim \int_{\infty}^x e^{\delta s}ds \lesssim e^{\delta x}, \quad \mbox{for} \quad x\leq 0,
\end{align*}
which gives  the result. The proof requires several tools of Fredholm theory for elliptic operators. The linearization of the equations \eqref{Main_equation} at $\Xi_{\infty}^{(c)}(\cdot,\cdot)$ gives
\begin{align}\label{Main_equation-Liouville_transform}
   \mathscr{L}_{\Xi}[v] :=  \Delta v + c\partial_x v + \left[\mu(x)  - 3 \left(\Xi_{\infty}^{(c)}(x,y)\right)^2\right]v,
\end{align}
 with domain of definition $ \mathcal{D}\left(\mathscr{L}_{\Xi}\right) = H^2\left(\mathbb{R}\times [0,\kappa]\right) \cap H_0^1\left(\mathbb{R}\times [0,\kappa]\right)$. Although this operator is nonself-adjoint, the limits as $|x|\to \infty$ of $\Xi^{(c)}$ are the same for all $c\in [0,2)$, therefore the results of \cite[Lem. 5.1 and 5.2]{Monteiro_Scheel} apply, showing that the operator $\widetilde{\mathscr{L}_{\Xi}}$ is Fredholm of index $0$, with essential spectrum strictly negative in $H^2\left(\mathbb{R}\times [0,\kappa]\right) \cap H_0^1\left(\mathbb{R}\times [0,\kappa]\right)$. According to \cite[Lem. 5.3]{Monteiro_Scheel}, we know that elements in the kernel of $\mathscr{L}_{\Xi}$ are spatially exponentially localized, namely, 
\begin{equation}\label{exponetial_decay:Kernel}
|\nabla u_0(x,y)|+|u_0(x,y )| \leq  C e^{- \delta|x|} ,\quad  (x, y)  \in  \mathbb{R}\times [0,\kappa] \quad a.e.\qquad \mbox{whenever} \quad u_0 \in \mathrm{Ker}\left(\mathscr{L}_{\Xi}\right)
\end{equation}
 Recall the partition of unity $\chi^{\pm}(\cdot)$ defined in \eqref{partition_of_unity}. We know that $\partial_x(\chi^{-}(x)\Xi^{(c)}(x,y)) \in \mathscr{C}^{\infty}(\R\times[0,\kappa];\R)$ solves a problem of the form 
$\displaystyle{\mathscr{L}_{\Xi}\left[v\right] = f},$
where $f\in L^2(\R\times[0,\kappa)$ is spatially localized. Writing $H^2(\R\times[0,\kappa]) = \mathrm{Ker}\left(\mathscr{L}_{\Xi}\right) \oplus \mathcal{X}$, we can assume with no loss of generality that $v\in \mathcal{X}$, thanks to property \eqref{exponetial_decay:Kernel} for elements in the kernel. However, as the operator $\mathscr{L}_{\Xi}: \mathcal{X} \to \mathrm{Rg}\left(\mathscr{L}_{\Xi}\right)$ is boundedly invertible, we can apply the same reasoning used in \cite[Cor. 5.5]{Monteiro_Scheel} to conclude that 
$e^{\delta<x>}v\in H^2(\R\times[0,\kappa]) \cap H_0^1(\R\times[0,\kappa])$ for all $\delta >0$ sufficiently small. A similar analysis can be done by considering $w(x,y) = \partial_y(\chi^{+}(y)\Xi^{(c)}(x,y))$, whence exponential rate of convergence to the far field as $y\to \infty$ is derived. It finishes the proof.
\end{Proof}

In fact, one can show by following the steps in the proof of  Lem. \ref{1_leadsto_0_1D:uniqueness} that the operator $\mathscr{L}_{\Xi}$ in \eqref{Main_equation-Liouville_transform} is boundedly invertible from $H^2(\R\times[0,\kappa])$ to  $L^2(\R\times[0,\kappa])$. Once more, using the IFT, we conclude the following result. 
 
\begin{Lemma}[Uniqueness of the  continuation in $c$;  $\mathcal{H}_{\kappa}$ problem, $\pi < \kappa < \infty$]\label{H_kappa:uniqueness_continuity_in_c} 
Recall the definition of $\mathcal{P}(c;\kappa)$ given in \eqref{critical_quantity}. For any fixed $\kappa \in (\pi, \infty)$ the following properties hold:
\begin{enumerate}
 \item there exists a unique solution $\Xi_{\kappa}^{(c)}(\cdot,\cdot)$ to the $\mathcal{H}_{\kappa}^{(c)}$ problem;
 \item the mapping   $c\mapsto \Xi_{\kappa}^{(c)}(\cdot):\left\{c\geq 0|\quad \mathcal{P}(c;\kappa) <1\right\} \to L^{\infty}(\R\times[0,\kappa];\R)$ is continuous.
\end{enumerate}
\end{Lemma}
\begin{Proof} The analysis is analogous to that of Lem. \ref{1_leadsto_0_1D:uniqueness} and is outline below, where we point out the necessary modifications. Fix $\kappa \in (\pi, \infty)$. Initially we define the linearized operator about the solutions $\Xi_{\kappa}^{c}$, obtaining the linearized operator
\begin{align*}
 \mathscr{L}_{\Xi^{(c)}}[v] = \partial_x^2 v + c\partial_x v + \mu(x) v  - 3(\Xi^{(c)})^2 v , \quad \mathcal{D}\left(\mathscr{L}_{\theta^{(c)}}\right) = H_2(\R\times[0,\kappa])\cap H_0^1(\R\times[0,\kappa]).
\end{align*}
Writing $v = e^{-\frac{c\,x}{2}} u$ we rewrite the above operator in a ``self-adjoint'' form,
\begin{align*}
 \widetilde{\mathscr{L}}_{\Xi^{(c)}}[v] = \partial_x^2 u + \left(\mu(x) - \frac{c^2}{4}\right) u  - 3(\Xi^{(c)})^2 u , \quad \mathcal{D}\left(\mathscr{L}_{\theta^{(c)}}\right) = e^{\frac{c\,x}{2}}\left(H_2(\R\times[0,\kappa])\cap H_0^1(\R\times[0,\kappa])\right).
\end{align*}
Notice that the mapping $u(\cdot) \mapsto e^{\frac{c\,x}{2}}u(\cdot)$ is an isometry between the spaces $H_2(\R\times[0,\kappa])\cap H_0^1(\R\times[0,\kappa])$ and $e^{\frac{c\,x}{2}}\left(H_2(\R\times[0,\kappa])\cap H_0^1(\R\times[0,\kappa])\right)$, therefore it suffices to  study the invertibility of $\widetilde{\mathscr{L}}_{\Xi^{(c)}}$ only. At this point we define $\widetilde{\mathscr{L}}_{\Xi^{(c)}}^{(\delta)}$ as the action of the operator $\widetilde{\mathscr{L}}_{\Xi^{(c)}}$ on the scale of Banach spaces $e^{\delta<x>}H^2(\R\times[0,\kappa])$. Continuity of the Fredholm index on the weight $\delta$ shows in the range $|\delta|\leq \frac{c}{2}$ we have that $\widetilde{\mathscr{L}}_{\Xi^{(c)}}^{(\delta)}$ are Fredholm operators with index zero, with the same kernel. As $\widetilde{\mathscr{L}}_{\Xi^{(c)}}^{(\delta)}\Big|_{\delta =0}$ is invertible one can argue as in Lem. \ref{1_leadsto_0_1D:uniqueness} to obtain the bounded invertibility of $\widetilde{\mathscr{L}}_{\Xi^{(c)}}$ in $e^{\frac{c\,x}{2}}\left(H_2(\R\times[0,\kappa])\cap H_0^1(\R\times[0,\kappa])\right).$ The results are then obtained from an application of the IFT.\end{Proof}

\subsection{Non-existence for the \texorpdfstring{$\mathcal{H}_{\kappa}$}{H-kappa}-problem: case \texorpdfstring{$\frac{c^2}{4} +\frac{\pi^2}{\kappa^2}>1$}{c2/4 + pi2/kappa2 >1} }
In this section we prove the non-existence of patterns when $\frac{c^2}{4} +\frac{\pi^2}{\kappa^2}>1$. The method is standard: 
roughly speaking, assuming an existing solution $\Xi_{\kappa}^{(c)}(\cdot,\cdot)$, we can obtain a supersolution $V(\cdot,\cdot)$ that is above $\Xi_{\kappa}^{(c)}(\cdot,\cdot)$ and which, under certain conditions, can touch the solution in at least one point. The function $z(\cdot,\cdot):= \Xi_{\kappa}^{(c)}(\cdot,\cdot) - V(\cdot,\cdot)$ solves an elliptic problem and has an interior maximum, which contradicts maximum principle and Hopf's Lemma. We make these words more precise in what follows. We begin with the main ingredients in the construction of the subsolution.

\begin{Lemma}[Non-existence; $\frac{c^2}{4} + \frac{\pi^2}{\kappa^2}>1$]
 Recall (see Fig. \ref{Figure:1d-solutions_no_jump}) the existence of a solution $w_d(\cdot)>0$ to the problem 
\begin{align} \label{equation:w_d}
\partial_x^2w_d + d\partial_xw_d + w_d - \left(w_d\right)^3, \qquad \lim_{x\to -\infty}w_d(x) = 1, \quad \lim_{x\to\infty}w_d(x) = 0,
\end{align}
where $d>\max\{2,c\}$ is chosen in such a way that  $\frac{d^2}{4} - \frac{c^2}{4} - \frac{\pi^2}{\kappa^2}<0$. Likewise, we write  $w_c(\cdot)$ for a function with similar spatial asymptotic properties and solving  $\partial_x^2w_c + c\partial_xw_c + w_c - \left(w_c\right)^3=0$. Define 
$$V(x,y) := w_c(x) + e^{-\frac{(c-d)\,x}{2}}w_d(x)v_{\kappa}(y), \qquad \mbox{where} \quad v_{\kappa}(y) := \alpha \sin\left(\frac{\pi\, y }{\kappa}\right).$$
Then $V(\cdot,\cdot)$ is a supersolution. Furthermore, $V(\cdot,\cdot)\leq \Xi_{\kappa}^{(c)}(\cdot,\cdot).$ 
\end{Lemma}
\begin{Proof}
Recall that $\mathcal{P}(c,\kappa) := \frac{c^2}{4} + \frac{\pi^2}{\kappa^2} >1$. The first step on our proof is the construction of a supersolution $V \geq \Xi_{\kappa}^{(c)}$. 
 Initially, we fix $w_c(\cdot)$ in such a way that $w_c(x) \geq 0$ on $x\leq 0$. Now, choosing $d>\max\{2,c\}$ in such a way that  $\frac{d^2}{4} - \frac{c^2}{4} - \frac{\pi^2}{\kappa^2}<0$ we fix a translated version of $w_d$ in such a way that 
\begin{align}\label{important_supersolution_inequ}
\frac{d^2}{4} - \frac{c^2}{4} - \frac{\pi^2}{\kappa^2} + w_d^2(x) - 3 w_c^2(x) \leq 0,\quad \mbox{whenever} \quad  x \leq 0. 
\end{align}
 In fact, we can exploit the translation invariance of solutions to the ODE \eqref{equation:w_d} and the fact that $x\mapsto w_{\{c,d\}}(x)$ is non-increasing for  $x\leq 0$ to conclude that  \eqref{important_supersolution_inequ} also holds for  translated version of $w_d(\cdot)$ and $w_c(\cdot)$ of the form $w_d(\cdot) \mapsto w_d(\cdot - \tau) $ and $w_c(\cdot) \mapsto w_c(\cdot - \tau)$, whenever $\tau \geq 0$. 
 
 Set $\widetilde{V}(x,y) = e^{-\frac{(c-d)\,x}{2}}w^{(d)}(x)v_{\kappa}(y)$ and write $V(x,y) = w_{c}(x) + \widetilde{V}(x,y)$. We claim that V is a supersolution on $x\leq 0$ (where $\mu(x) \equiv 1$). Indeed, a tedious but straightforward computation shows that
\begin{align*}
\Delta_{x,y}  V + c\partial_xV + V - V^3 &=  \widetilde{V}\left[\frac{d^2}{4} - \frac{c^2}{4} - \frac{\pi^2}{\kappa^2} \right] + e^{-\frac{(c-d)\,x}{2}}w_{d}^3(x)v_{\kappa}(y) + w_c^3 - \left[w_c(x) + e^{-\frac{(c-d)\,x}{2}}w_d(x)v_{\kappa}(y)\right]^3\\
& =   \widetilde{V}\left[\frac{d^2}{4} - \frac{c^2}{4} - \frac{\pi^2}{\kappa^2} +  w_d^2(x) - 3 w_c^3(x)\right] + J_2
\end{align*}
where $J_2= J_2(w_c, w_d, v_{\kappa})$ consists of non-positive terms only. Thus,
\begin{align}\label{V_is_supersolution}
\Delta_{x,y}  V + c\partial_xV + \mu(x)V - V^3 &=   \widetilde{V}\left[\frac{d^2}{4} - \frac{c^2}{4} - \frac{\pi^2}{\kappa^2} +  w_d^2 - 3 w_c^3\right] + J_2 \leq 0, \quad \mbox{in}\quad x\leq 0,
\end{align}
thanks to \eqref{important_supersolution_inequ}, which concludes the proof that V is a supersolution in $x\leq 0$.

Our second step consists in proving that $V \geq \Xi^{(c)}$. Unfortunately, we cannot play with the translation of both $w_d$ and $w_c$ separately without destroying the inequality in \eqref{V_is_supersolution}, so we pursue another route. The asymptotic properties of $w_c$ imply that $V(x,y) \geq w_c(x) \geq \Xi^{(c)}(x,y)$ as $x \to -\infty$. As $w_c(\cdot)$ can be chosen so that $w_c(0) = 0$ we need to understand the term $\widetilde{V}(x,y)$ for $x <0$, $x$ finite. Inspecting the calculations that lead to \eqref{V_is_supersolution}, one can see that the parameter $\alpha \geq 0$ in  $v_{\kappa}(y):= \alpha \sin\left(\frac{\pi\,y}{\kappa}\right)$ plays no crucial role. Standard elliptic regularity estimates show that $\partial_y\Xi^{(c)}(\cdot,\cdot)$ is a  bounded function on $\R\times[0,\kappa]$. Consequently, we can take $\alpha>0$  sufficiently large so that 
\begin{align}\label{estimate_supersolution_on_x_equal_0}
 w_d(0) v_{\kappa}(y) = \alpha w_d(0) \sin\left(\frac{\pi \, y}{\kappa}\right)\geq \Xi^{(c)}(0,y);
\end{align}
thanks to  $w_d(\cdot) >0$. Furthermore, the inequality \eqref{estimate_supersolution_on_x_equal_0} persists  when we use a translated version of $w_d(\cdot)$, for $x \mapsto  w_d(x)$ is non-increasing, namely, $ w_d(-\tau) v_{\kappa}(y) \geq w_d(0) v_{\kappa}(y)\geq \Xi^{(c)}(0,y),$ for all $\tau >0$. Thus, for a fixed $\tilde{M}>0$ we have that 
$$\widetilde{V}(x,y) \geq \Xi^{(c)}(x,y),\quad \mbox{whenever} \quad x\in [-\tilde{M},0]\times [0,\kappa].$$
Translating $w_c(\cdot)$ in such a way that $w_c(x) \geq \Xi^{(c)}(x,y)$ for $x \leq -\tilde{M}$ and $y \in [0,\kappa]$ we get that $V(x,y)\geq \Xi^{(c)}(x,y)$ for all $(x,y)\in (-\infty,0]\times[0,\kappa]$. The third and last part consists of varying the parameter   $\alpha >0$ continuously in order to satisfy $V \geq \Xi^{(c)}$ with an equality in at least one point. Now the proof goes as in Prop. \eqref{Prop:Dichotomy}: define $Z(x,y) := \Xi(x,y) - \bar{V}(x,y)$; clearly $Z(x,y)\leq 0$ in $x\leq 0$, $y \in [0,\kappa]$.   The function $Z$ satisfies
\begin{align*}
 \Delta_{x,y} Z + c\partial_x Z + \mu(x)Z - f[\Xi,\bar{V}]\left(Z\right) \geq 0, \quad \mbox{on} \quad x\leq 0.
\end{align*}
where $f[a,b]:= \frac{a^3 - b^3}{a-b}$ whenever $a\neq b$ and $3a^2$ otherwise.  We can apply classical maximum principles, since for $M>0$ chosen large enough in order to give $M + f[\Xi,\bar{V}] >0$; thus we write the above equation as 
\begin{align*}
 \Delta_{x,y} Z + c\partial_x Z + \mu(x)Z - (f[\Xi,\bar{V}] + M)\left(Z\right)\geq - M Z \geq 0
\end{align*}
As $Z$ has a maximum point in $\R\times[0,\kappa]$ we get that $Z\equiv 0$. This contradiction leads to the non-existence of solutions in the case $\frac{c^2}{4} +\frac{\pi^2}{\kappa^2}>1$.   
\end{Proof}

Finally, we put all these auxiliary results together and prove the main result of this section:
\begin{Proof}[of Theorem \ref{theorem:2D_result}; case $\pi <\kappa < \infty$] Combine the above discussion of the non-existence of solutions with the results of Lemmas \ref{H_kappa:existence}, \ref{H_kappa:exponential_convergence}, and  \ref{H_kappa:uniqueness_continuity_in_c}
\end{Proof}

\section{Two-dimensional quenched patterns -- single horizontal interfaces; \texorpdfstring{$\mathcal{H}_{\infty}$}{H infty} problem}\label{Sec:2D_result:k_infinite}
In this section, we shall prove Theorem \ref{theorem:2D_result} in the case $\kappa =\infty$. To be consistent with the notation introduced in Section \ref{Sec:H_1d_and_1_0_problem:truncation_approach} we exploit the fact that the nonlinearity in \eqref{Main_equation} is odd to solve the problem in the half space $\R\times (-\infty,0]$. Further symmetries of the equation are also exploited: we solve an equivalent $\mathcal{H}_{\infty}$-problem, seeking for a solution  $\Xi_{\infty}^{(c)}(\cdot,\cdot)$ to \eqref{Main_equation} in $\R^2$, satisfying
\begin{equation}\label{H_infty:asymptotic_limits}
\lim_{x\to -\infty}\Xi_{\infty}^{(c)}(x,y)= -\tanh\left(\frac{y}{\sqrt{2}}\right),\quad \lim_{y\to 
\pm\infty}\Xi_{\infty}^{(c)}(x,y)=\mp\theta^{(c)}(x),\quad \mbox{and}\quad \lim_{x\to \infty}\Xi_{\infty}^{(c)}(x,y)=0,
\end{equation}
where $\theta^{(c)}(\cdot)$ is the one-dimensional solution to the $\onelc$-problem. Notice that the odd symmetry solutions to \eqref{Main_equation} with respect to  $\left(y, \Xi_{\infty}^{(c)}(x,y)\right)\mapsto \left(-y,-\Xi_{\infty}^{(c)}(x,-y)\right)$ readily gives the pattern with the properties stated in Theorem \ref{theorem:2D_result}. Moreover, the Dirichlet boundary conditions at $y=0$ removes the non-uniqueness of solutions induced by $y$-translation invariance. 

 Thanks to the results of Sec. \ref{Sec:H_1d_and_1_0_problem:truncation_approach} related to the 1-d problem the following observation is readily available. 
\begin{Observation}[Restriction to the case $c<2$]
 It is clear from \eqref{H_infty:asymptotic_limits} that the above problem is meaningless when $c>2$ for the patterns $\theta^{(c)}$ do not exist. We can readily say that no solution to this problem exists when $c>2$, immediately restricting our study to the range $0 \leq c<2$.
 \end{Observation}

In fact, in this section we prove that for all quenching fronts speeds in the range $c\in[0,2)$ there exists a unique $\mathcal{H}_{\infty}$-pattern (up to translations in the $y$ direction), which corresponds to the statement of Theorem\ref{theorem:2D_result}. The strategy goes as in Sections \ref{Sec:H_1d_and_1_0_problem:truncation_approach} and  \ref{Sec:2D_result:k_finite}: first, by reducing the problem to a half plane and truncating it, restricting the problem to a rectangle $\Omega_{(-M,L)} := (-M,L)\times(-M, 0)$. Then we let $M\to\infty$ and, subsequently, we let $L\to\infty$. For the sake of simplicity in this section we will omit any sub-index $\infty$.

The \emph{truncated} $\mathcal{H}_\infty$-problem is set up as 
\begin{align*}
\left\{\begin{array}{rclr}
	  \Delta_{x,y} u +c\partial_xu + \mu(x) u - u^3&=&0, & (x,y)\in\Omega_{(-M,L)},\nonumber\\
u&=&g_{(-M,L)},& (x,y)\in\partial\Omega_{(-M,L)},\label{2D_H_infty_problem}
       \end{array}\right.
\end{align*}
where $g_{(-M,L)}(x,y):=\theta_{(-M,L)}^{(c)}(x)\cdot \theta_{(-M,0)}^{(c)}(y)$,   $\theta_{(-M,L)}^{(c)}(\cdot)$ the solutions to the truncated one-dimensional problem  \eqref{H_1D_problem:truncated} on the interval $(-M, L)$. Similarly to the results of Lem. \ref{H_kappa:existence_uniqueness} and \cite[Sec. 3]{Monteiro_Scheel}, we construct unique solutions to this truncated problem using an iterative scheme. The solution $\Theta_{(-M,L)}^{(c)}(\cdot,\cdot): \Omega_{(-M,L)} \to [0,1]$ is shown to be unique; Furthermore, exploiting that $0\leq u\leq 1$ and  Agmon-Douglis-Nirenberg regularity we readily conclude that $u$ and derivatives are H\"older continuous across $x=0$, and, in fact, $\Theta_{(-M,L)}^{(c)}\in \mathscr{C}^{(1,\alpha)}(\overline{\Omega_{(-M,L)}}), $  for all $0\leq\alpha <1$.
Following the method in Sec. \ref{Sec:2D_result:k_finite}, we extend these functions to the whole plane $\R^2$:
\begin{equation*}\label{extension_operator:H_infty}
 \mathscr{E}\left[\Theta_{(-M,L)}^{(c)}\right](x,y) = \left \{\begin{array}{cl}
                          \Theta_{(-M,L)}^{(c)}(x,y), &   \quad  \mbox{for} \quad  (x,y) \in \Omega_{(-M,L)},\\
                           \mathscr{E}\left[\theta_{(-M,L)}^{(c)}\right](x)\cdot \mathscr{E}\left[\theta_{(-M,0)}^{(c)}\right](y),&\quad  \mbox{for} \quad (x,y) \in \mathbb{R}^2 \setminus \Omega_{(-M,L)}.
                          \end{array}\right.   \end{equation*} 
We summarize the main properties of the functions $\mathscr{E}\left[\Theta_{(-M,L)}^{(c)}\right](\cdot,\cdot) $ in the following Proposition, whose proof is similar to that of \cite[Prop. 3.4]{Monteiro_Scheel}:
\begin{Proposition}[Properties of the extension operator, $\mathcal{H}_{\infty}$-problem]\label{Thm:monotonicity_properties_H_infty_problem} The following properties of  $\mathscr{E}[\cdot]$ hold.
\begin{enumerate}[label=(\roman*), ref=\theTheorem(\roman*)]
\hangindent\leftmargin
 \item \emph{(Monotonicity of $\mathscr{E}$)}\label{Thm:monotonicity_properties_H_infty_problem:monotonicity_extension_operator}
  We have $0 \leq \mathscr{E}\left[\Theta_{(-M,L)}^{(c)}\right](\cdot, \cdot) \leq 1$. Furthermore, if $w$ is only defined in a subset $A \subset \mathbb{R}^2$ so that  $ \Omega_{(-M,L)} \subset A \subset \Omega_{(-\infty,L)}$, $0 \leq w(\cdot) \leq 1$ and $w(\cdot) \leq \Theta_{(-M,L)}^{(c)}(\cdot, \cdot)$,  then 
$0 \leq \mathscr{E}\left[w\right](\cdot) \leq \mathscr{E}\left[\Theta_{(-M,L)}^{(c)}(\cdot, \cdot)\right](\cdot)  \, \, \mbox{in} \, \, \mathbb{R}^2 $.                                                                                                                                                                                                                                        
 \item \emph{(Monotonicity in $M$)}\label{Thm:monotonicity_properties_H_infty_problem:2D_monotonicity_in_M}
Let  $0 \leq M <\widetilde{M}$ and $L \geq 0$ be fixed. Then $
M <\widetilde{M} \implies   \mathscr{E}\left[U_{(-\widetilde{M},L)}\right](x,y) \leq  \mathscr{E}\left[\Theta_{(-M,L)}^{(c)}\right](x,y). $
 \item \emph{(Monotonicity in $L$)} \label{Thm:monotonicity_properties_H_infty_problem:2D_monotonicity_in_L}
Let  $0 \leq L <\widetilde{L}$ and $M \geq 0$ be fixed. Then $
L <  \widetilde{L} \implies   \mathscr{E}\left[\Theta_{(-M,L)}^{(c)}\right](x,y) \leq  \mathscr{E}\left[U_{(-M,\widetilde{L})}\right](x,y). $
 \item \emph{(Monotonicity in $x$)}\label{Thm:monotonicity_properties_H_infty_problem:2D_monotonicity_in_x}
Let  $L,M,y$be fixed. Then
the mapping  $ x \to \mathscr{E}\left[\Theta_{(-M,L)}^{(c)}\right](x,y)$ is non-increasing. 
\item \emph{(Monotonicity in $y$)}\label{Thm:monotonicity_properties_H_infty_problem:2D_monotonicity_in_y}
Let  $L,M,x$be fixed. Then
the mapping  $ y \to \mathscr{E}\left[\Theta_{(-M,L)}^{(c)}\right](x,y)$ is non-increasing. 
\end{enumerate}
\end{Proposition}

\begin{Lemma}[Existence; $\mathcal{H}_{\infty}$ problem]\label{H_infty:existence} There exists a solution $\Xi_{\infty}^{(c)}(\cdot, \cdot)$ satisfying  \eqref{Main_equation} in $\R^2$ satisfying the limits \eqref{H_infty:asymptotic_limits}. Furthermore, $\partial_y \Xi_{\infty}^{(c)}(\cdot, \cdot) \leq 0.$ 
\end{Lemma}
\begin{Proof}
We obtain a solution in $(x,y) \in \R\times(-\infty,0)$, extending it to the whole plane $\R^2$ as a function satisfying $\Xi^{(c)}(x,-y) = -\Xi^{(c)}(x,y)$. The result follows the construction in Lem. \ref{H_kappa:existence}, first defining $\displaystyle{\Theta_{(-\infty,L)}(\cdot,\cdot) :=\inf_{M>0}\Theta_{(-\infty,L)}(\cdot,\cdot)}$ and then defining  
 \begin{align*}
  \Xi_{\infty}^{(c)}(\cdot,\cdot) = \sup_{L>0}\Theta_{(-\infty,L)}(\cdot,\cdot).
 \end{align*}
 Further properties of $\Xi_{\infty}^{(c)}(\cdot, \cdot)$ are derived as in Lem. \ref{1_leadsto_0_1D:existence} (see also \cite[\S 5]{Monteiro_Scheel}). The inequality $\partial_y \Xi_{\infty}^{(c)}(\cdot, \cdot) \leq 0$ is a consequence of  Prop. \ref{Thm:monotonicity_properties_H_infty_problem:2D_monotonicity_in_y}, because the mapping $y \mapsto \Xi_{\infty}^{(c)}(x, y)$ is non-increasing.  
\end{Proof}

The next result provides a crucial ingredient in the construction of interfaces with contact angle:
\begin{Corollary}[Monotonicity in the y-direction]\label{Monotonicity_in_y}
 For all $x\in \R$ the mapping $\displaystyle{y \mapsto  \Xi_{\infty}^{(c)}(x,y) }$ is strictly monotonic. 
\end{Corollary}
\begin{Proof}
 It suffices to show that $\partial_y\Xi_{\infty}^{(c)}(x,\cdot)\neq 0$ a.e., thanks to the regularity of $\Xi_{\infty}^{(c)}(\cdot,\cdot)$. The inequality $\partial_y\Xi_{\infty}^{(c)} \leq 0$ is obtained without difficulty by making use of Prop. \ref{Thm:monotonicity_properties_H_infty_problem:2D_monotonicity_in_y} and the limiting construction in Lem. \ref{H_infty:existence}. The strict inequality $\partial_y\Xi_{\infty}^{(c)} < 0$ a.e. then follows using a Harnack inequality (cf. \cite[Theorem 9.22]{gilbarg2015elliptic}; see also \cite[Prop. 4.1]{Monteiro_Scheel-contact_angle}).  
\end{Proof}
\begin{Lemma}[Exponential convergence; $\mathcal{H}_{\infty}$ problem]\label{H_infty:exponential_convergence} The limits in \eqref{H_infty:asymptotic_limits} take place at exponential rate, i.e., there exists a $C$, $\delta >0$ independent of $x$ and $y$ such that
\begin{align*}
\lim_{x\to \infty}|\Xi_{\infty}^{(c)}(x,y)| \leq Ce^{-\delta |x|}, \quad  \lim_{x\to \infty}\left|\Xi_{\infty}^{(c)}(x,y) - \tanh\left(\frac{y}{\sqrt{2}}\right)\right| \leq Ce^{-\delta |x|} ,\mbox{and} \lim_{y\to \pm\infty}|\Xi_{\infty}^{(c)}(x,y) \pm \theta^{(c)}(x)|\leq Ce^{-\delta |y|}.
\end{align*}
\end{Lemma}
\begin{Proof}
 The proof is similar to that of Lem. \ref{H_kappa:exponential_convergence}, with small differences. As before, we first trace the essential spectrum of  the linearized operator \eqref{Main_equation} at $\Xi_{\infty}^{(c)}(\cdot,\cdot)$ gives
\begin{align}\label{H_infty:Main_equation-Liouville_transform}
   \mathscr{L}_{\Xi}[v] :=  \Delta v + c\partial_x v + \left[\mu(x)  - 3 \left(\Xi_{\infty}^{(c)}(x,y)\right)^2\right]v,
\end{align}
 with domain of definition $ \mathcal{D}\left(\mathscr{L}_{\Xi}\right) = H_{\mathrm{odd}}^2\left(\mathbb{R}^2\right):= \{w \in H^2\left(\mathbb{R}^2\right); w(x,y) = -w(x,-y)\}.$ We remark that \cite[Lem. 5.3]{Monteiro_Scheel} still applies: indeed, we can define the asymptotic operators
   \begin{align*}
\mathscr{M}^-[v] := \partial_y^2v(\cdot)  + \left[1- 3\tanh^2\left(\frac{y}{\sqrt{2}}\right)\right]v(\cdot),  \qquad  \mathscr{M}^+[v] :=  \partial_y^2v(\cdot)  -v(\cdot),        
 \end{align*}
with domain  $\mathcal{D}\left(\mathscr{M}^{\pm} \right) =  H_{\mathrm{odd}}^2([0,\kappa])\cap H_0^1([0,\kappa]).$ We claim that these operators are invertible: indeed,  coercivity implies that $\mathrm{Ker}\left(\mathscr{M}^{-}\right) = \{0\}$; The same holds in the case of $\mathscr{M}^{+}$, for in $H^2$ its kernel is given by  $\{\partial_y\tanh\left(\frac{\cdot}{\sqrt{2}}\right)\}$ which is a simple eigenfunction, since it has no nodal points (see \cite[Theorem 8.38]{gilbarg2015elliptic}).  Therefore, in the space of odd functions, we have that both operators $\mathscr{M}^{\pm}$ are invertible.  In our next step, we argue as in \cite[Lem. 5.1]{Monteiro_Scheel}, describing limiting operators associated with $\mathscr{L}_{\Xi}[\cdot]$:
\begin{subequations}
\begin{align}
 \mathscr{L}_{\Xi}^{(x\to +\infty)}[v] &= \left(\partial_x^2 + c\partial_x + \mathscr{M}^+\right)[v], \label{H_infty:far_field:x_minus_inf}\\
 \mathscr{L}_{\Xi}^{(x\to -\infty)}[v] &= \left(\partial_x^2 + c\partial_x + \mathscr{M}^-\right)[v],\label{H_infty:far_field:x_plus_inf}\\
 \mathscr{L}_{\Xi}^{(y\to +\infty)}[v]&= \left(\Delta_{x,y} + c\partial_x + \mu(x) - 3(\theta^{(c)})^2\right)[v]\label{H_infty:far_field:y_plus_inf}. 
\end{align}
\end{subequations}
Fourier transforming the operators \eqref{H_infty:far_field:x_minus_inf}-\eqref{H_infty:far_field:x_plus_inf} in $x$, and \eqref{H_infty:far_field:y_plus_inf} in $y$, shows that these operators are boundedly invertible. Observe that $v(x,y) = \partial_x\left(\chi^-(x)\Xi_{\infty}^{(c)}(x,y)\right)$ (resp., $v(x,y) = \partial_x\left(\chi^+(x)\Xi_{\infty}^{(c)}(x,y)\right)$) satisfies  $\mathscr{L}_{\Xi}^{(x\to -\infty)}[v] = f$ (resp., $\mathscr{L}_{\Xi}^{(x\to +\infty)}[v] = f$) where the left hand side is spatially localized, it follows from the same reasoning as in Lem. \ref{H_kappa:exponential_convergence} that $e^{\delta<x>}v \in H^2(\R^2)\cap H_0^1(\R^2)$. A similar result holds for 
 $v = \partial_y (\chi^+(y)\Xi)$, which solves  $\mathscr{L}_{\Xi}^{(y\to +\infty)}[v] = f$ for $f$ spatially localized in $y$, i.e., $f\in e^{\delta <y>}L^2(\R^2)$, where $<y> = \sqrt{1 +y^2}$. The rest of the analysis is similar to that used in the proof of Lem. \ref{H_kappa:exponential_convergence}.
\end{Proof}
\begin{Lemma}[Uniqueness up to translation in the y-direction]\label{H_infty:uniqueness} Whenever $c\geq 0$ the solutions constructed in Theorem \ref{theorem:2D_result} and the solutions constructed in \cite{Monteiro_Scheel} by continuation are the same up to translation in the $y$-direction.\end{Lemma}
\begin{Proof}
 The machinery given in \cite{Monteiro_Scheel-contact_angle} can be used  to derive a simple proof: first notice that  $\partial_y \Xi_{\infty}\in Ker\mathscr{L}_{\Xi_{\infty}}[\cdot]$. As $\partial_y\Xi$ has a sign one can use \cite[Lem. 4.9]{Monteiro_Scheel-contact_angle} to conclude that $\partial_y\Xi^{\infty}=C\partial_y\Theta $ for a $C$ constant. Upon integration in y and using the fact that both solutions converge to the same limit as $y \to \pm \infty$ and satisfy $\Theta(x,0) = \Xi^{\infty}(x,0) = 0$ we get that $C \equiv 1$, and the result follows. 
\end{Proof}

Our last result concerns the continuity of the patterns $\Xi_{\infty}^{(c)}(\cdot,\cdot)$ in $c$.

\begin{Proposition}\label{H_infty:continuity_in_c}
 The mapping $\displaystyle{c\mapsto \Xi_{\infty}^{(c)}(\cdot):[0,2) \to L^{\infty}(\R^2;\R)}$ is continuous.
\end{Proposition}
\begin{Proof}
As before, we exploit the symmetries of the problem to reduce the analysis to the half-space $\R\times(-\infty,0]$. A similar analysis to that of Lem. \ref{1_leadsto_0_1D:uniqueness} shows that the linearization of the equation \eqref{Main_equation} at $\Xi_{\infty}^{(c)}$ given in \eqref{H_infty:Main_equation-Liouville_transform} is an invertible operator from $H^2(\R\times(-\infty,0])\cap H_0^1(\R\times(-\infty,0])$ (i.e., restricted to odd solutions in $y$) to $L^2(\R\times(-\infty,0])$. However, the result does not readily follows, because
$$ \Xi_{\infty}^{(c)} - \Xi_{\infty}^{(d)} \not \in H^2(\R\times(-\infty,0])\cap H_0^1(\R\times(-\infty,0]).$$
 We circumvent this issue with a far-field core decomposition, which we now explain: for a fixed $d\in[0,2)$,  we use an Ansatz of the form  $ w(x,y)+ \Xi_{\infty}^{(c)} + \chi^-(y) \left[\theta^{(d)}(x) - \theta^{(c)}(x)\right]$, where $w \in H^{2}(\R\times(-\infty,0];\R)$; notice that $\Xi_{\infty}^{(c)} + \chi^-(y) \left[\theta^{(d)}(x) - \theta^{(c)}(x)\right]$ converges exponentially fast to its far-field and asymptotically solves the PDE \eqref{Main_equation}. We rewrite problem \eqref{Main_equation} as
 \begin{align*}
  \mathscr{L}_{\Xi_{\infty}}[w] : = \Delta_{x,y} w + c\partial_x w + \mu(x) w - 3(\Xi_{\infty}^{(c)})^2 w = \mathscr{N}\left(w, \Xi_{\infty}^{(c)},\theta^{(d)}(x), \theta^{(c)}(x) \right)
 \end{align*}
A reasoning similar to that in the proof of Lem. \ref{H_kappa:uniqueness_continuity_in_c} shows that $\mathscr{L}_{\Xi_{\infty}}$ is a bounded, invertible operator from $H^2(\R\times(-\infty,0])\cap H_0^1(\R\times(-\infty,0])$ to $L^2(\R\times(-\infty,0])$. Furthermore, the right hand side is in $L^2(\R\times(-\infty,0])$, thanks to the Sobolev embedding $H^2(\R^2) \hookrightarrow L^{\infty}(\R^2)$. The conclusion then follows from the IFT.
\end{Proof}

 \begin{Proof}[of Theorem \ref{theorem:2D_result}; case $\kappa = \infty$] Combine the results of Lemmas \ref{H_infty:existence}, \ref{H_infty:exponential_convergence}, \ref{H_infty:uniqueness}, and Prop. \ref{H_infty:continuity_in_c}.
 \end{Proof}  
\section{Discussion}\label{Discussion} Among several possible directions of further investigation, we would like to mention the following:

\paragraph{Metastability of patterns.} As addressed by the numerical studies of  \cite[Sec. VI]{Foard}, defining the  parameter regions of metastability for creation of patterns (either perpendicular or parallel to the quenching front) is a challenging and interesting direction of investigation. From a broader perspective, a numerical, if possible analytical,  description of parameter curves on the boundary of different morphological states would be valuable in applications.

\paragraph{Selection mechanisms.} What are the crucial mechanisms involved in the wavenumber selection in the wake of the front? How relevant are the nonlinearity and the speed of the quenching front in this selection?
We refer to \cite[\S 3.3]{Nishiura} for a general discussion about wavenumber selection.

\paragraph{Critical cases; \texorpdfstring{$\mathcal{P}(c;\kappa) = 1$}{S=1}.} The behavior of the patterns in this critical scenario possibly requires a different approach, since one can see in the proofs of Theorem \ref{theorem:1D_result} and \ref{theorem:2D_result} that the speed of the quenching front has to be away from the critical case. The result would be interesting and add valuable knowledge in the classification of patterns obtained from directional quenching.

\paragraph{Non-planar quenching fronts and oblique stripes.} 
Is it possible to control the contact angle of the $\hp_{\kappa}$-patterns? Although it was shown in \cite{Monteiro_Scheel} that oblique patterns do not exist in \eqref{Main_equation}, these patterns can still exist in the case of the unbalanced equation \eqref{unbalanced_equation}. It is worth to mention that  the result of \cite{Monteiro_Scheel-contact_angle} describes a family of solutions displaying one single (almost) horizontal interface whose  contact angle with the quenching front can be varied by modification of the chemical potential parameters across the quenching interface. We refer to   \cite{Racz} for physics motivation and a more detailed discussion on the chirality of helicoidal patterns in the context of recurrent precipitation. 
%

\newcommand{\etalchar}[1]{$^{#1}$}

\end{document}